\documentclass[aap,seceqn,citesort,MSNbibl,dvips]{arximspdf}

\doi{10.1214/10-AAP677}
\volume{20}
\issue{6}
\pubyear{2010}
\firstpage{2086}
\lastpage{2117}

\makeatletter

\newcommand{\bd}{\bolds}

\newproclaim{definition}{Definition}[section]
\newtheorem{theorem}{Theorem}
\newtheorem{lemma}{Lemma}[section]
\newtheorem{coro}{Corollary}[section]
\newproclaim{remark}{Remark}

\makeatother

\begin{document}
\begin{frontmatter}

\title{Spectral distributions of adjacency and Laplacian matrices
of random graphs}
\runtitle{Spectral of Laplacian matrices}

\begin{aug}
\author[A]{\fnms{Xue} \snm{Ding}\ead[label=e1]{dingxue@stat.umn.edu}} and
\author[B]{\fnms{Tiefeng} \snm{Jiang}\corref{}\thanksref{t1}\ead[label=e2]{tjiang@stat.umn.edu}}
\runauthor{X. Ding and T. Jiang}
\affiliation{Jilin University and University of Minnesota, and
University of Minnesota}
\address[A]{School of Mathematics\\
Jilin University\\
2699 Qianjin Street\\
Changchun\\
China\\
and\\
School of Statistics\\
University of Minnesota\\
224 Church Street\\
Minneapolis, Minnesota 55455\\
USA\\
\printead{e1}} 
\address[B]{School of Statistics\\
University of Minnesota\\
224 Church Street\\
Minneapolis, Minnesota 55455\\
USA\\
\printead{e2}}
\end{aug}

\thankstext{t1}{Supported in part by NSF Grant DMS-04-49365.}

\received{\smonth{4} \syear{2009}}
\revised{\smonth{11} \syear{2009}}

%
\begin{abstract}
In this paper, we investigate the spectral properties of the adjacency
and the Laplacian matrices of random graphs. We prove that:

\hspace*{5.21pt}(i) the law of large numbers for the spectral norms
and the
largest eigenvalues of the adjacency and the Laplacian matrices;

\hspace*{2.86pt}(ii) under some further independent conditions, the
normalized largest eigenvalues of the Laplacian matrices are dense in a
compact interval almost surely;

(iii) the empirical distributions of the eigenvalues of the Laplacian
matrices converge weakly to the free convolution of the standard
Gaussian distribution and the Wigner's semi-circular law;

\hspace*{1pt}(iv) the empirical distributions of the eigenvalues of the
adjacency matrices converge weakly to the Wigner's semi-circular law.
\end{abstract}

%

\end{frontmatter}

\section{Introduction}\label{intro}

The theory of random graphs was founded in the late 1950s by Erd\"{o}s
and R\'{e}nyi
\cite{ER59,ER60,ER61a,ER61b}.
The work of Watts and Strogatz \cite{WAST} and Barab\'{a}si and Albert
\cite{BAAL} at the end of the last century initiated new interest in
this field.
The subject is at the intersection between graph theory and probability
theory. One can see, for example,
\cite{Bollobas,Chung97,CL06,Collind,Durrett06,ES74,JLR,Kol,palmer} for
book-length treatments.

The spectral graph theory is the study of the properties of a graph in
relationship to the characteristic polynomial, eigenvalues and
eigenvectors of its adjacency matrix or Laplacian matrix. For
reference, one can see books
\cite{Chung97,Puppe} for the deterministic case and \cite{CL06} for
the random case, and literatures therein. The spectral graph theory has
applications in chemistry
\cite{BLW} where eigenvalues were relevant to the stability of
molecules. Also, graph spectra appear naturally in numerous questions
in theoretical physics and quantum mechanics (see, e.g.,
\cite{Evan92,EvanEc,Evan83,MF91,Novikov,RB,RDD}). For connections
between the eigenvalues of the adjacency matrices and the Laplacian
matrices of graphs and Cheeger constants, diameter bounds, paths and
routing in graphs, one can see
\cite{CL06}.

Although there are many matrices
for a given graph with $n$ vertices, the most studied are their
adjacency matrices and the Laplacian matrices. Typically,
random graphs are considered with the number of vertices $n$ tending
to infinity. Many geometrical and topological properties can be
deduced for a large class of random graph ensembles, but the spectral
properties of the random graphs are still uncovered to a large extent.

In this paper, we will investigate the spectral properties of the
adjacency and the Laplacian matrices of some random graphs.
The framework of the two matrices will be given next.

Let $n\geq2$ and $\Gamma_{n}=(\mathcal{V}_{n},E_{n})$ be a graph, where
$\mathcal{V}_{n}$ denotes a set of $n$ vertices $v_{1}, v_{2},
\ldots, v_{n}$, and $E_{n}$ is the set of edges. In this paper, we
assume that the edges in $E_n$ are always nonoriented. For basic
definitions of graphs, one can see, for example,
\cite{Bollo}. The adjacency matrix and the Laplacian matrix of the
graph are of the form
%
%
\begin{equation}\label{Dec}
\mathbf{A}_{n}=\pmatrix{
0 & \xi_{12}^{(n)}& \xi_{13}^{(n)}&\cdots& \xi_{1n}^{(n)}\vspace*{1pt}\cr
\xi_{21}^{(n)}& 0 & \xi_{23}^{(n)}&\cdots& \xi_{2n}^{(n)}\vspace*{1pt}\cr
\xi_{31}^{(n)}& \xi_{32}^{(n)}& 0 &\cdots& \xi_{3n}^{(n)}\vspace*{1pt}\cr
\vdots&\vdots&\vdots&\vdots&\vdots\vspace*{1pt}\cr
\xi_{n1}^{(n)}& \xi_{n2}^{(n)}& \xi_{n3}^{(n)}&\cdots& 0}
\end{equation}
and
%
%
\begin{equation}\label{Mar}
\bd{\Delta}_{n}=\pmatrix{
\displaystyle\sum_{j\neq1}\xi_{1j}^{(n)}&-\xi_{12}^{(n)}&-\xi_{13}^{(n)}&\cdots
&-\xi_{1n}^{(n)}\vspace*{2pt}\cr
-\xi_{21}^{(n)}&\displaystyle\sum_{j\neq2}\xi_{2j}^{(n)}&-\xi_{23}^{(n)}&\cdots
&-\xi_{2n}^{(n)}\vspace*{2pt}\cr
-\xi_{31}^{(n)}&-\xi_{32}^{(n)}&\displaystyle\sum_{j\neq3}\xi_{3j}^{(n)}&\cdots
&-\xi_{3n}^{(n)}\vspace*{2pt}\cr
\vdots&\vdots&\vdots&\vdots&\vdots\cr
-\xi_{n1}^{(n)}&-\xi_{n2}^{(n)}&-\xi_{n3}^{(n)}&\cdots&\displaystyle\sum_{j\neq
n}\xi_{nj}^{(n)}}
\end{equation}
with relationship
%
%
\begin{equation}\label{gameb}
\bd{\Delta}_n=\mathbf{D}_n-\mathbf{A}_n,
\end{equation}
where $\mathbf{D}_n=(\sum_{l\neq i}^{n}\xi_{il}^{(n)})_{1\leq i \leq n}$
is a diagonal matrix.

As mentioned earlier, we will focus on nonoriented random graphs in
this paper. Thus, the adjacency matrix $\mathbf{A}_{n}$ is always
symmetric. If the graph is also simple, the entry $\xi_{ij}^{(n)}$ for
$i\ne j$ only takes value $1$ or $0$ with $1$ for an edge between $v_i$
and~$v_j$, and $0$ for no edge between them.

The Laplacian matrix $\bd{\Delta}_n$ for graph $\Gamma_n$ is also
called the admittance matrix or the Kirchhoff matrix in literature. If
$\Gamma_n$ is a simple random graph, the $(i, i)$-entry of $\bd
{\Delta}_n$ represents the degree of vertex $v_i$, that is, the number
of vertices connected to $v_i$. $\bd{\Delta}_n$ is always
nonnegative\vadjust{\goodbreak} (this is also true for $\bd{\Delta}_n$ as long as the
entries $\{\xi_{ij}^{(n)}; 1\leq i\ne j \leq n\}$ are nonnegative);
the smallest eigenvalue of $\bd{\Delta}_n$ is zero; the second
smallest eigenvalue stands for the algebraic connectivity; the
Kirchhoff theorem establishes the relationship between the number of
spanning trees of $\Gamma_n$ and the eigenvalues of $\bd{\Delta}_n$.

An \textit{Erd\"{o}s--R\'{e}nyi random graph} $G(n, p)$ has $n$ vertices.
For each pair of vertices $v_i$ and $v_j$ with $i\ne j$, an edge
between them is formed randomly with chance $p_n$ and independently of
other edges (see
\cite{ER59,ER60,ER61a,ER61b}). This random graph corresponds to
Bernoulli entries $\{\xi_{ij}^{(n)}; 1\leq i< j\leq n\}$, which are
independent random variables with $P(\xi_{ij}^{(n)}=1)=1-P(\xi
_{ij}^{(n)}=0)=p_n$ for all $1\leq i< j \leq n$.

For weighted random graphs, $\{\xi_{ij}^{(n)}; 1\leq i< j\leq n\}$
are independent random variables and $\xi_{ij}^{(n)}$ is a product of
a Bernoulli random variable $\operatorname{Ber}(p_n)$ and a nice random variable, for
instance, a Gaussian random variable or a random variable with all
finite moments (see, e.g.,
\cite{KSV,KKPS}). For the sign model studied in \cite
{BG01,KKPS,RB,RDD}, $\xi_{ij}^{(n)}$ are independent random variables
taking three values: $0, 1, -1$. In this paper, we will study the
spectral properties of $\mathbf{A}_n$ and $\bd{\Delta}_n$ under more
general conditions on $\{\xi_{ij}^{(n)}; 1\leq i< j\leq n\}$ [see
(\ref{shufic})].

Now we need to introduce some notation about the eigenvalues of
matrices. Given an $n\times n$ symmetric matrix $\mathbf{M}$. Let $
\lambda_{1}\geq\lambda_{2}\geq\cdots\geq\lambda_{n}$ be the
eigenvalues of~$\mathbf{M}$, we sometimes also write this as $ \lambda
_{1}(\mathbf{M})\geq\lambda_{2}(\mathbf{M})\geq\cdots\geq\lambda
_{n}(\mathbf{M})$ for clarity. The notation
$\lambda_{\max}=\lambda_{\max}(\mathbf{M})$, $\lambda_{\min}=\lambda
_{\min}(\mathbf{M})$ and $\lambda_{k}(\mathbf{M})$ stand for the largest
eigenvalue, the smallest eigenvalue and the $k$th largest eigenvalue of
$\mathbf{M}$, respectively. Set
%
%
\begin{eqnarray}\label{minyi}
\hat{\mu}(\mathbf{M})&=&\frac{1}{n}\sum_{i=1}^n\delta_{\lambda_i}
\quad\mbox{and}\nonumber\\[-8pt]\\[-8pt]
F^{\mathbf{M}}(x)&=&\frac{1}{n}\sum_{i=1}^{n}I(\lambda
_{i}\leq x),\qquad x \in\mathbb{R}.\nonumber
\end{eqnarray}
Then, $\hat{\mu}(\mathbf{M})$ and $F^{\mathbf{M}}(x)$ are the empirical spectral
distribution of $\mathbf{M}$ and the empirical spectral cumulative
distribution function of $\mathbf{M}$, respectively.


In this paper, we study $\mathbf{A}_n$ and $\bd{\Delta}_n$ not only for
random graphs but also study them in the context of random matrices.
Therefore, we allow the entries $\xi_{ij}^{(n)}$'s to take real values
and possibly with mean zero. It will be clear in our theorems if the
framework is in the context of random graphs or that of of random matrices.

Under general conditions on $\{\xi_{ij}^{(n)}\}$, we prove in this
paper that a suitably normalized $\hat{\mu}(\mathbf{A}_{n})$ converges
to the semi-circle law; a suitably normalized $\hat{\mu}(\bd{\Delta
}_{n})$ converges weakly to the free convolution of the standard normal
distribution and the semi-circle law. Besides, the law of large numbers
for largest eigenvalues and the spectral norms of $\mathbf{A}_{n}$ and
$\bd{\Delta}_{n}$ are obtained. Before stating these results, we need
to give the assumptions on the entries of $\mathbf{A}_n$ in (\ref{Dec})
and $\bd{\Delta}_n$ in (\ref{Mar}).
%
%
\begin{eqnarray}\label{shufic}
\begin{tabular}{p{315pt}}
Let $\{\xi_{ij}^{(n)};1\leq i< j\leq n, n\geq2\}$
be random variables defined on the same probability space and
$\{\xi_{ij}^{(n)};1\leq i< j\leq n\}$ be independent for
each $n\geq2$ (not necessarily identically distributed)
with $\xi_{ij}^{(n)}=\xi_{ji}^{(n)}$, $E(\xi
_{ij}^{(n)})=\mu_n$, $\operatorname{Var}(\xi_{ij}^{(n)})=\sigma_n^2>0$
for all $1\leq i< j\leq n$ and $n\geq2$ and
$\sup_{1\leq i<j\leq n, n\geq2}E|(\xi_{ij}^{(n)}-\mu_n)/\sigma
_n|^{p}<\infty$
for some $p>0$.
\end{tabular}\hspace*{-48pt}
\end{eqnarray}
The values of $p$ above will be specified in each result later. In what
follows, for an $n\times n$ matrix $\mathbf{M}$, let $\|\mathbf{M}\|=\sup
_{\mathbf{x}\in\mathbb{R}^n\dvtx\|\mathbf{x}\|=1}\|\mathbf{M}\mathbf{x}\|$ be the
spectral norm of $\mathbf{M}$, where $\|\mathbf{x}\|=\sqrt{x_1^2+\cdots+
x_n^2}$ for $\mathbf{x}=(x_1, \ldots, x_n)' \in\mathbb{R}^n$. Now we
state the main results of this paper.
\begin{theorem}\label{west1} Suppose (\ref{shufic}) holds for some
$p>6$. Assume $\mu_n=0$ and $\sigma_n=1$ for all $n\geq2$. Then:
%
\begin{enumerate}[(a)]
\item[(a)] $\frac{\lambda_{\max}(\bd{\Delta}_{n})}{\sqrt{n\log n}}\to
\sqrt{2}$ in probability as $n\to\infty$.
\end{enumerate}
Furthermore, if $\{\bd{\Delta}_2, \bd{\Delta}_3, \ldots\}$ are
independent, then:
\begin{enumerate}[(a)]
\item[(b)]
$\liminf_{n\to\infty}\frac{\lambda_{\max}(\bd
{\Delta}_{n})}{\sqrt{n\log n}}=\sqrt{2}$ a.s. and
$\limsup_{n\to\infty}\frac{\lambda_{\max}(\bd{\Delta
}_{n})}{\sqrt{n\log n}}=2$ a.s., and
the sequence $\{\lambda_{\max}(\bd{\Delta}_{n})/\sqrt
{n\log n}; n\geq2\}$ is dense in $[\sqrt{2}, 2]$
a.s.;
\item[(c)] the conclusions in \textup{(a)} and \textup{(b)} still hold if
$\lambda_{\max}(\bd{\Delta}_{n})$ is replaced by $\|\bd
{\Delta}_{n}\|$.
\end{enumerate}
\end{theorem}

For typically-studied random matrices such as the Hermite ensembles and
the Laguerre ensembles, if we assume the sequence of $n\times n$
matrices for all $n\geq1$ are independent as in Theorem \ref{west1},
the conclusions (b) and (c) in Theorem \ref{west1} do not hold. In
fact, for Gaussian Unitary Ensemble (GUE), which is a special case of
the Hermite ensemble, there is a large deviation inequality
$P(|n^{-1/2}\lambda_{\max}-\sqrt{2}|\geq\varepsilon)\leq
e^{-nC_{\varepsilon}}$ for any $\varepsilon>0$ as $n$ is sufficiently
large, where $C_{\varepsilon}>0$ is some constant (see (1.24) and
(1.25) from
\cite{Ledoux} or \cite{BADG}). With or without the independence
assumption, this inequality implies from the Borel--Cantelli lemma that
$n^{-1/2}\lambda_{\max}\to\sqrt{2}$ a.s. as $n\to\infty$.
Similar large deviation inequalities also hold for Wishart and sample
covariance matrices (see, e.g., \cite{FHJ,VMB}).

For two sequence of real numbers $\{a_n; n\geq1\}$ and $\{b_n; n\geq
1\}$, we write $a_n\ll b_n$ if $a_n/b_n\to0$ as $n\to\infty$, and
$a_n\gg b_n$ if $a_n/b_n\to+\infty$ as $n\to\infty$. We use $n\gg
1$ to denote that $n$ is sufficiently large.
\begin{coro}\label{coro1} Suppose (\ref{shufic}) holds for some $p>6$.
Then, as $n\to\infty$:
%
\begin{enumerate}[(a3)]
\item[(a1)] $\frac{\lambda_{\max}(\bd{\Delta}_{n})}{\sigma_n\sqrt
{n\log n}}\to\sqrt{2}$ in probability if
$|\mu_n|\ll\sigma_n (\frac{\log n}{n} )^{1/2}$;

\item[(a2)] $\frac{\lambda_{\max}(\bd{\Delta}_{n})}{n\mu_n}\to1$
in probability if $\mu_n>0$ for $n\gg1$ and
$\mu_n\gg\sigma_n (\frac{\log n}{n} )^{1/2}$;

\item[(a3)] $\frac{\lambda_{\max}(\bd{\Delta}_{n})}{n\mu_n}\to0$
in probability if $\mu_n<0$ for $n\gg1$
and $|\mu_n|\gg\sigma_n (\frac{\log n}{n} )^{1/2}$.
\end{enumerate}
Furthermore, assume $\{\bd{\Delta}_2, \bd{\Delta}_3, \ldots\}$
are independent, then:
\begin{enumerate}[(b3)]
\item[(b1)] $\liminf_{n\to\infty}\frac{\lambda_{\max}(\bd
{\Delta}_{n})}{\sigma_n\sqrt{n\log n}}=\sqrt{2}$ a.s.
and $\limsup_{n\to\infty}\frac{\lambda_{\max}(\bd{\Delta
}_{n})}{\sigma_n\sqrt{n\log n}}=2$ a.s., and the sequence
$\{\frac{\lambda_{\max}(\bd{\Delta}_{n})}{\sigma_n\sqrt
{n\log n}}; n\geq2 \}$ is dense in $[\sqrt{2}, 2]$
a.s. if $|\mu_n|\ll\sigma_n (\frac{\log n}{n}
)^{1/2}$;
\item[(b2)] $\lim_{n\to\infty}\frac{\lambda_{\max}(\bd
{\Delta}_{n})}{n\mu_n}= 1$ a.s. if $\mu_n>0$ for
$n\gg1$ and $\mu_n\gg\sigma_n (\frac{\log n}{n} )^{1/2}$;
\item[(b3)] $\lim_{n\to\infty}\frac{\lambda_{\max}(\bd{\Delta
}_{n})}{n\mu_n}= 0$ a.s. if $\mu_n<0$ for $n\gg1$
and $|\mu_n|\gg\sigma_n (\frac{\log n}{n} )^{1/2}$.
\end{enumerate}
Finally, \textup{(a1)} and \textup{(b1)} still hold if $\lambda_{\max}(\bd{\Delta
}_{n})$ is replaced by $\|\bd{\Delta}_n\|$; if $\xi_{ij}^{(n)}\geq
0$ for all $i, j, n$, then \textup{(a2)} and \textup{(b2)} still hold if $\lambda_{\max
}(\bd{\Delta}_{n})$ is replaced by $\|\bd{\Delta}_n\|$.
\end{coro}
\begin{remark}\label{Remark1}
For the Erd\"{o}s--R\'{e}nyi random
graph, the condition ``(\ref{shufic}) holds for some $p>p_0$'' with
$p_0>2$ is true only when $p_n$ is bounded away from zero and one. So,
under this condition of $p_n$, Corollary \ref{coro1} holds. Moreover,
under the same restriction of $p_n$, Theorems \ref{free1} and \ref
{Ruth}, that will be given next, also hold.


Let $\{\nu, \nu_1, \nu_2, \ldots\}$ be a sequence of probability
measures on $\mathbb{R}$. We say that $\nu_n$ converges weakly to
$\nu$ if $\int_{\mathbb{R}}f(x) \nu_n(dx)\to\int_{\mathbb
{R}}f(x) \nu(dx)$ for any bounded and continuous function $f(x)$
defined on $\mathbb{R}$. The Portmanteau lemma says that the weak
convergence can also be characterized in terms of open sets or closed
sets (see, e.g., \cite{Dudley}).


Now we consider the empirical distribution of the eigenvalues of the
Laplacian matrix $\bd{\Delta}_n$. Bauer and Golinelli \cite{BG01}
simulate the eigenvalues for the Erd\"{o}s--R\'{e}nyi random graph with
$p$ fixed. They observe that the limit $\nu$ of the empirical
distribution of $\lambda_i(\bd{\Delta}_n), 1\leq i \leq n$, has a
shape between the Gaussian and the semicircular curves. Further, they
conclude from their simulations that $m_4/m_2^2$ is between $2$ and
$3$, where $m_i$ is the $i$th moment of probability measure $\nu$. In
fact, we have the following result.
\end{remark}
\begin{theorem}\label{free1} Suppose (\ref{shufic}) holds for some
$p>4$. Set $\tilde{F}_n(x)=\frac{1}{n}\times\break\sum_{i=1}^nI \{\frac
{\lambda_i(\bd{\Delta}_n)-n\mu_n}{\sqrt{n}\sigma_n}\leq x \}$
for $x\in\mathbb{R}$.
Then, as $n\to\infty$, with probability one, $\tilde{F}_n$ converges weakly
to the free convolution $\gamma_{M}$ of the
semicircular law and the standard normal distribution. The measure
$\gamma_{M}$
is a nonrandom symmetric probability measure with
smooth bounded density, does not depend on the distribution of
$\{\xi_{ij}^{(n)}; 1\leq i<j\leq n, n\geq2\}$ and has an unbounded
support.
\end{theorem}


More information on $\gamma_{M}$ can be found in \cite{BDJ}. For the
Erd\"{o}s--R\'{e}nyi random graphs, the weighted random graphs in
\cite{KSV,KKPS} and the sign models in \cite{BG01,KKPS,RB,RDD}, if
$p_n$ is bounded away from $0$ and $1$ as $n$ is large, then (\ref
{shufic}) holds for all $p>4;$ thus Theorem \ref{free1} holds for all
of these graphs.

It is interesting to notice that the limiting curve appeared in Theorem
\ref{free1} is indeed a hybrid between the standard Gaussian
distribution and the semi-circular law, as observed in
\cite{BG01}. Moreover, for the limiting distribution, it is shown in
\cite{BDJ} that $m_4/m_2^2=8/3\in(2,3)$, which is also consistent
with the numerical result in
\cite{BG01}.

Before introducing the next theorem, we now make a remark. It is proved
in \cite{BDJ} that the conclusion in the above theorem holds when $\xi
_{ij}^{(n)}=\xi_{ij}$ for all $1\leq i< j \leq n$ and $n\geq2$, where
$\{\xi_{ij}; 1\leq i<j<\infty\}$ are independent and identically
distributed random variables with $E\xi_{12}=0$ and $E(\xi
_{12})^2=1$. The difference is that the current theorem holds for any
independent, but not necessarily identically distributed, random
variables with arbitrary mean $\mu_n$ and variance $\sigma_n^2$.


Now we consider the adjacency matrices. Recall $\mathbf{A}_n$ in (\ref
{Dec}). Wigner \cite{Wigner} establishes the celebrated semi-circle
law for matrix $\mathbf{A}_n$ with entries $\{\xi_{ij}^{(n)}=\xi_{ij}\dvtx
1\leq i< j <\infty\}$ being i.i.d. $N(0,1)$-distributed random
variables (for its extensions, one can see, e.g.,
\cite{BS} and literatures therein). Arnold \cite{Arnold67,Arnold71}
proves that Wigner's result holds also for the entries being i.i.d.
random variables with a finite sixth moment. In particular, this
implies that, for the adjacency matrix $\mathbf{A}_n$ of the
Erd\"{o}s--R\'{e}nyi random graph with $p$ fixed, the empirical
distribution of the
eigenvalues of $\mathbf{A}_n$ converges to the semi-circle law (see also Bollobas
\cite{Bollobas}). In the next result we show that, under a condition
slightly stronger than a finite second moment, the semicircular law
still holds for $\mathbf{A}_n$.

\begin{theorem}\label{Becky}
Let $\omega_{ij}^{(n)}:=(\xi_{ij}^{(n)}-\mu_n)/\sigma_n$ for all
$i,j,n$. Assume (\ref{shufic}) with $p=2$ and
\[
\max_{1\leq i<j\leq
n}E \bigl\{\bigl(\omega_{ij}^{(n)}\bigr)^{2}I \bigl(\bigl|\omega_{ij}^{(n)}\bigr|\geq
\varepsilon
\sqrt{n} \bigr) \bigr\}\to0
\]
as $n\to\infty$ for any $\varepsilon>0$, which is particularly true
when (\ref{shufic}) holds for some $p>2$. Set
\[
\tilde{F}_n(x)=\frac{1}{n}\sum_{i=1}^nI \biggl\{\frac{\lambda_i(
\mathbf{A}_n) + \mu_n}{\sqrt{n}\sigma_n}\leq x \biggr\},\qquad x\in\mathbb{R}.
\]
Then, almost surely, $\tilde{F}_n$ converges weakly to the
semicircular law with
density $\frac{1}{2\pi}\sqrt{4-x^{2}} I(|x|\leq2)$.
\end{theorem}


Applying Theorem \ref{Becky} to the Erd\"{o}s--R\'{e}nyi random graph,
we have the following result.
%
\begin{coro}\label{cheng} Assume (\ref{shufic}) with $P(\xi
_{ij}^{(n)}=1)=p_n=1-P(\xi_{ij}^{(n)}=0)$ for all $1\leq i<j \leq n$
and $n\geq2$. If $\alpha_n:=(np_n(1-p_n))^{1/2}\to\infty$ as $n\to
\infty$, then, almost surely, $F^{\mathbf{A}_n/\alpha_n}$ converges
weakly to the semicircular law with
density $\frac{1}{2\pi}\sqrt{4-x^{2}}I(|x|\leq2)$. In particular,
if $1/n\ll p_n\to0$ as $n\to\infty$, then, almost surely, $F^{
\mathbf{A}_n/\sqrt{np_n}}$ converges weakly to the same semicircular law.
\end{coro}
%

The condition ``$\alpha_n:=(np_n(1-p_n))^{1/2}\to\infty$ as $n\to
\infty$'' cannot be relaxed to that ``$np_n\to\infty$.'' This is
because, as $p_n$ is very close to $1$, say, $p_n=1$, then $\xi
_{ij}^{(n)}=1$ for all $i\ne j$. Thus $\mathbf{A}_n$ has eigenvalue $n-1$
with one fold and $-1$ with $n-1$ fold. This implies that $F^{
\mathbf{A}_n} \to\delta_{-1}$ weakly as $n\to\infty$.

Corollary \ref{cheng} shows that the semicircular law holds not only
for $p$ being a constant as in Arnold \cite{Arnold67,Arnold71}, it
also holds for the dilute Erd\"{o}s--R\'{e}nyi graph, that is, $1/n\ll
p_n\to0$ as $n\to\infty$. A result in Rogers and Bray
\cite{RB} (see also a discussion for it in\vspace*{-2pt} Khorunzhy
et al. \cite{KKPS}) says that, if $P(\xi
_{ij}^{(n)}=\pm1)=p_n/2$ and $P(\xi_{ij}^{(n)}=0)=1-p_n$, the
semicircular law holds for the corresponding $\mathbf{A}_n$ with $1/n\ll
p_n\to0$. It is easy to check that their result is a corollary of
Theorem \ref{Becky}.

Now we study the spectral norms and the largest eigenvalues of
$\mathbf{A}_n$. For the Erd\"{o}s--R\'{e}nyi random graph, the largest
eigenvalue of $\mathbf{A}_n$ is studied in
\cite{FKO,KrSu}. In particular, following Juh\'{a}z \cite{Juh},
F\"{u}redi and Koml\'{o} \cite{FKO} showed that the largest eigenvalue
has asymptotically a normal distribution when $p_n=p$ is a constant;
Krivelevich and Sudakov
\cite{KrSu} proved a weak law of large numbers for the largest
eigenvalue for the full range of $p_n\in(0,1)$. In the following, we
give a result for $\mathbf{A}_n$ whose entries do not necessarily take
values of $0$ or $1$ only.
Recall $\lambda_k(\mathbf{A}_n)$ and $\|\mathbf{A}_n\|$ are the $k$th
largest eigenvalue and the spectral norm of $\mathbf{A}_n$, respectively.
%
\begin{theorem}\label{Ruth}
Assume (\ref{shufic}) holds for some $p>6$. Let $\{k_n; n\geq1\}$
be a sequence of positive integers such that $k_n=o(n)$ as $n\to\infty$.
The following hold:

\begin{longlist}
\item If $\lim_{n\to\infty}\mu_n/(n^{-1/2}\sigma_{n})= 0$, then $\|
\mathbf{A}_n\|/\sqrt{n}\sigma_n \to2$ a.s. and $\lambda_{k_n}
(\mathbf{A}_n)/\break(\sqrt{n}\sigma_n)\to2$ a.s. as $n\to\infty$.

\item If $\lim_{n\to\infty}\mu_{n}/(n^{-1/2}\sigma_{n})=+\infty$,
then $\lambda_{\max}(\mathbf{A}_{n})/(n\mu_{n})\to1$ a.s. as $n\to
\infty$.

\item If $\lim_{n\to\infty}|\mu_{n}|/(n^{-1/2}\sigma_{n})=+\infty$,
then $\|\mathbf{A}_{n}\|/(n|\mu_{n}|)\to1$ a.s. as $n\to\infty$.
\end{longlist}
\end{theorem}
%
%
\begin{remark}\label{Remark2} The conclusion in (ii) cannot be
improved in general to that $\lambda_{k_n}(\mathbf{A}_{n})/(n\mu_{n})\to
1$ a.s. as $n\to\infty$. This is because when $\sigma_n$ is
extremely small, $\mathbf{A}_n$ roughly looks like $\mu_n(\mathbf{J}_n -
\mathbf{I}_n)$, where all the entries of $\mathbf{J}_n$ are equal to one, and
$\mathbf{I}_n$ is the $n\times n$ identity matrix. It is easy to see that
the largest eigenvalue of $\mu_n(\mathbf{J}_n-\mathbf{I}_n)$ is $(n-1)\mu
_n>0$, and all of the remaining $n-1$ eigenvalues are identical to $-1$.


From the above results, we see two probability distributions related to
the spectral properties of the random graphs: the Wigner's semi-circle
law and the free convolution of the standard normal distribution and
the semi-circle law. The Kesten--McKay law is another one. It is the
limit of the empirical distributions of the \mbox{eigenvalues} of the random
$d$-regular graphs (see
\cite{McKay}).

The proofs of Theorems \ref{west1} and \ref{free1} rely on the moment
method and some tricks developed in \cite{BDJ}. Theorems \ref{Becky}
and \ref{Ruth} are derived through a general result from
\cite{bai99} and certain truncation techniques in probability theory.

The rest of the paper is organized as follows: we will prove the
theorems stated above in the next section; several auxiliary results
for the proofs are collected in the \hyperref[appendix]{Appendix}.
\end{remark}

\section{Proofs}\label{proofs}

%
\begin{lemma}\label{Mikea}
Let $\mathbf{U}_n=(u_{ij}^{(n)})$ be an $n\times n$ symmetric random
matrix, and $\{u_{ij}^{(n)}; 1\leq i\leq j \leq n, n\geq1\}$ are
defined on the same probability space. Suppose, for each $n\geq1$, $\{
u_{ij}^{(n)};1\leq i\leq j\leq n\}$ are independent random variables
with $Eu_{ij}^{(n)}=0, \operatorname{Var}(u_{ij}^{(n)})=1$
for all\vspace*{1pt} $1\leq i,j\leq n$, and $\sup_{1\leq i,j\leq
n, n\geq1}E|u_{ij}^{(n)}|^{6+\delta}<\infty$ for some $\delta>0$. Then:
%
\begin{longlist}
\item $\lim_{n\to\infty}\frac{\lambda_{\max}(\mathbf{U}_n)}{\sqrt
{n}}=2$
a.s. and
$\lim_{n\to\infty}\frac{\|\mathbf{U}_n\|}{\sqrt{n}}=2$ a.s.;
\item
the statements in \textup{(i)} still hold if $\mathbf{U}_n$
is replaced by $\mathbf{U}_n - \operatorname{diag}(u_{ii}^{(n)})_{1\leq i
\leq n}$.
\end{longlist}
\end{lemma}

The proof of this lemma is a combination of Lemmas \ref{missing} and
\ref{sprc} in the \hyperref[appendix]{Appendix} and some truncation techniques. It is
postponed and will be given later in this section.
\begin{pf*}{Proof of Theorem \ref{west1}}
First, assume (a) and (b) hold. Since $\mu_n=0$ for all $n\geq2$, (a)
and (b) also hold if $\lambda_{\max}(\bd{\Delta}_n)$ is replaced by
$\lambda_{\max}(-\bd{\Delta}_n)$. From the symmetry of $\bd{\Delta}_n$,
we know that
\[
\|\bd{\Delta}_n\|=\max\{-\lambda_{\min}(\bd{\Delta}_n), \lambda
_{\max}(\bd{\Delta}_n)\}=\max\{\lambda_{\max}(-\bd{\Delta}_n),
\lambda_{\max}(\bd{\Delta}_n)\}.
\]
Now the function $h(x, y):=\max\{x, y\}$ is continuous in $(x, y)\in
\mathbb{R}^2$, applying the two assertions
\[
\limsup_{n\to\infty}\max\{a_n, b_n\}=\max \Bigl\{
\limsup_{n\to\infty}a_n, \limsup_{n\to\infty}b_n \Bigr\}
\]
and
\[
\liminf_{n\to\infty}\max\{a_n, b_n\}\geq\max
\Bigl\{\liminf_{n\to\infty}a_n, \liminf_{n\to\infty}b_n \Bigr\}
\]
for any $\{a_n\in\mathbb{R}; n\geq1\}$ and $\{b_n\in\mathbb
{R}; n\geq1\}$, we obtain $\|\bd{\Delta}_n\|/\sqrt{n\log n}$
converges to $\sqrt{2}$ in probability, and
\[
\liminf_{n\to\infty}\frac{\|\bd{\Delta}_{n}\|}{\sqrt{n\log
n}}\geq\sqrt{2}\qquad\mbox{a.s.} \quad\mbox{and}\quad \limsup_{n\to\infty
}\frac{\|\bd{\Delta}_{n}\|}{\sqrt{n\log n}}=2 \qquad\mbox{a.s.}
\]
and
\[
\mbox{the sequence } \biggl\{\frac{\|\bd{\Delta}_{n}\|}{\sqrt
{n\log n}}; n\geq2 \biggr\} \mbox{ is dense in } \bigl[\sqrt{2}, 2\bigr]\qquad
\mbox{a.s.}
\]
Thus (c) is proved. Now we turn to prove (a) and (b).

Recall (\ref{gameb}), $\bd{\Delta}_n=\mathbf{D}_n-\mathbf{A}_n$.
First, $\lambda_{\max}(\mathbf{D}_{n})-\|\mathbf{A}_{n}\| \leq\lambda
_{\max}(\bd{\Delta}_{n})\leq\break
\lambda_{\max}(\mathbf{D}_{n})+\|\mathbf{A}_{n}\|$ for all $n\geq2$. Second,
by (ii) of Lemma \ref{Mikea}, $\|\mathbf{A}_n\|/\sqrt{n}\to2$ a.s.
as $n\to\infty$.
Thus, to prove (a) and (b) in the theorem, it is enough to show that
%
%
\begin{eqnarray}\qquad
\label{probasha}
& \displaystyle\frac{T_n}{\sqrt{n\log n}}\to\sqrt{2} \qquad\mbox{in
probability};&\\
\label{mada}
& \displaystyle\liminf_{n\to\infty}\frac{T_n}{\sqrt{n\log n}}=\sqrt{2}
\qquad\mbox{a.s.} \quad\mbox{and}\quad \limsup_{n\to\infty}\frac{T_n}{\sqrt{n\log
n}}=2\qquad\mbox{a.s.;}&
\\
\label{africa}
& \displaystyle\mbox{the sequence } \biggl\{\frac{T_n}{\sqrt{n\log n}}; n\geq
2 \biggr\} \mbox{ is dense in } \bigl[\sqrt{2}, 2\bigr] \qquad\mbox{a.s.,}&
\end{eqnarray}
where $T_n=\lambda_{\max}(\mathbf{D}_n)=\max_{1\leq i \leq n}\sum_{j\ne
i}\xi_{ij}^{(n)}$ for $n\geq2$.\vspace*{10pt}

\textsc{Proof of (\ref{probasha})}.\quad By Lemma \ref{mamam}, for each
$1\leq i \leq n$ and $n\geq2$, there exist i.i.d. $N(0,1)$-distributed
random variables $\{\eta_{ij}^{(n)}; 1\leq j\leq n, j \ne i\}$ for
each $n\geq2$ such that
%
%
\begin{eqnarray}\label{tus}
&&
\max_{1\leq i\leq n}P \biggl(\biggl|\sum_{j\ne i}\xi_{ij}^{(n)}-\sum_{ j\ne i}
\eta_{ij}^{(n)}\biggr|\geq\varepsilon\sqrt{n\log n} \biggr) \nonumber\\
&&\qquad\leq
\frac{C}{1+(\varepsilon\sqrt{n\log
n})^{6}}\sum_{j\ne i}E\bigl|\xi_{ij}^{(n)}\bigr|^{6}\\
&&\qquad\leq\frac{C}{n^{2}(\log n)^{3}},\nonumber
\end{eqnarray}
where here and later in all proofs, $C$ stands for a constant not
depending on $i,j$ or~$n$, and may be different from line to line. It is
well known that
%
%
\begin{equation}\label{zhengt}
\frac{x}{\sqrt{2\pi}(1+x^{2})}e^{-x^{2}/2}\leq P\bigl(N(0,1)\geq x\bigr)\leq
\frac{1}{\sqrt{2\pi}x}e^{-x^{2}/2}
\end{equation}
for any $x>0$. Since $\sum_{j\ne i}\xi_{ij}^{(n)}\leq\sum_{j\ne
i}\eta_{ij}^{(n)} + |\sum_{j\ne i}\xi_{ij}^{(n)}- \sum_{j\ne i}\eta
_{ij}^{(n)}|$, then
%
%
\begin{eqnarray}
&&P\bigl(T_n\geq(\alpha+ 2\varepsilon)\sqrt{n\log n}\bigr)\nonumber\\
&&\qquad
\leq n\cdot\max_{1\leq i \leq n}P\biggl(\sum_{j\ne i}\xi
_{ij}^{(n)}\geq(\alpha+ 2\varepsilon)\sqrt{n\log n}
\biggr)\nonumber\\[-8pt]\\[-8pt]
&&\qquad \leq n\cdot\max_{1\leq i \leq n}P\biggl(\sum_{j\ne i}\eta
_{ij}^{(n)}\geq(\alpha+ \varepsilon)\sqrt{n\log n} \biggr)\nonumber\\
&&\qquad\quad{} + n\cdot\max_{1\leq i \leq n}P\biggl(\biggl|\sum_{j\ne i}\xi_{ij}^{(n)}-
\sum_{j\ne i}\eta_{ij}^{(n)}\biggr|\geq\varepsilon\sqrt{n\log n} \biggr)\nonumber
\end{eqnarray}
for any $\alpha>0$ and $\varepsilon>0$. Noticing $\sum_{j\ne i}\eta
_{ij}^{(n)}\sim\sqrt{n-1}\cdot N(0,1)$ for any $1\leq i \leq n$, by
(\ref{zhengt}) and then (\ref{tus}),
%
%
\begin{eqnarray}\label{lanse}\qquad
P\bigl(T_n\geq(\alpha+ 2\varepsilon)\sqrt{n\log n}\bigr)
& \leq& n P\bigl(N(0,1)\geq(\alpha+ \varepsilon)\sqrt{\log n} \bigr) + \frac
{C}{n(\log n)^3}\nonumber\\[-8pt]\\[-8pt]
& \leq& C n^{1-(\alpha+\varepsilon)^2/2} + \frac{C}{n(\log
n)^3}\nonumber
\end{eqnarray}
for $n$ sufficiently large. In particular, taking $\alpha=\sqrt{2}$,
we obtain that
%
%
\begin{equation}\label{dedao}
P \biggl(\frac{T_n}{\sqrt{n\log n}} \geq\sqrt{2}+ 2\varepsilon
\biggr)=O \biggl(\frac{1}{n^{\varepsilon}} \biggr)
\end{equation}
as $n\to\infty$ for any $\varepsilon\in(0,1]$, since the last term in
(\ref{lanse}) is of order $n^{-1}(\log n)^{-3}$ as $n\to\infty$.

Define $k_{n}=[n/\log n]$ and $V_{n}={\max_{1\leq i\leq k_{n}}}|\sum
_{1\leq j\leq k_{n}
}\xi_{ij}^{(n)}|$ with $\xi_{ii}^{(n)}=0$ for all $1\leq i \leq n$.
By the same argument as in obtaining (\ref{lanse}), we have that, for
any fixed $\alpha>0$,
%
%
\begin{equation}\label{houh}
P \biggl(\frac{V_n}{\sqrt{k_n\log k_n}} \geq\alpha+ 2\varepsilon
\biggr)\leq C (k_n)^{1-(\alpha+\varepsilon)^2/2} + \frac{C}{k_n(\log k_n)^3}
\end{equation}
as $n$ is sufficiently large. Noticing $n/k_n\to\infty$, and taking
$\alpha+\varepsilon=10$ above, we have
%
%
\begin{equation}\label{shengchu}
P\bigl(V_n\geq\varepsilon\sqrt{n\log n}\bigr)\leq\frac{1}{n(\log n)^{3/2}}
\end{equation}
as $n$ is sufficiently large. Observe that
%
%
\begin{equation}\label{most}
T_{n}\geq\max_{1\leq i\leq
k_{n}}\sum_{j=k_{n}+1}^{n}\xi_{ij}^{(n)}-V_{n}.
\end{equation}
Similarly to (\ref{tus}), by Lemma \ref{mamam}, for each $1\leq i
\leq n$ and $n\geq2$, there exist i.i.d. $N(0,1)$-distributed random
variables $\{\zeta_{ij}^{(n)}; 1\leq i\leq n, j\ne i\}$ such that
%
%
\begin{eqnarray}\label{mue}
&&
\max_{1\leq i\leq k_n}P \Biggl(\Biggl|\sum_{j=k_{n}+1}^{n}\xi
_{ij}^{(n)}-\sum_{j=k_{n}+1}^{n}
\zeta_{ij}^{(n)}\Biggr|\geq\varepsilon\sqrt{n\log n} \Biggr) \nonumber\\
&&\qquad\leq
\frac{C}{1+(\varepsilon\sqrt{n\log
n})^{6}}\sum_{j=k_{n}+1}^{n}E\bigl|\xi_{ij}^{(n)}\bigr|^{6}\\
&&\qquad\leq\frac{C}{n^{2}(\log n)^{3}}\nonumber
\end{eqnarray}
as $n$ is sufficiently large for any $\varepsilon>0$. Fix $\beta>0$. By
(\ref{most}), (\ref{shengchu}) and then independence
%
%
\begin{eqnarray}\label{meishu}\qquad
& & P\bigl(T_n\leq(\beta-2\varepsilon)\sqrt{n\log n}\bigr)\nonumber\\
&&\qquad\leq P\Biggl(\max_{1\leq i\leq
k_{n}}\sum_{j=k_{n}+1}^{n}\xi_{ij}^{(n)}\leq(\beta-\varepsilon)\sqrt
{n\log n}\Biggr) + P\bigl(V_n\geq\varepsilon\sqrt{n\log n}\bigr)\\
&&\qquad\leq \max_{1\leq i\leq
k_{n}}P\Biggl(\sum_{j=k_{n}+1}^{n}\xi_{ij}^{(n)}\leq(\beta-\varepsilon
)\sqrt{n\log n}\Biggr)^{k_n} + \frac{1}{n(\log n)^{3/2}}\nonumber
\end{eqnarray}
as $n$ is sufficiently large. Observe that
%
%
\begin{eqnarray}\label{bede}
& & P\Biggl(\sum_{j=k_{n}+1}^{n}\xi_{ij}^{(n)}\leq(\beta-\varepsilon)\sqrt
{n\log n}\Biggr)\nonumber\\
&&\qquad \leq P\Biggl(\sum_{j=k_{n}+1}^{n}\zeta_{ij}^{(n)}\leq\biggl(\beta-\frac
{\varepsilon}{2}\biggr)\sqrt{n\log n}\Biggr)\\
&&\qquad\quad{} + P\Biggl(\Biggl|\sum_{j=k_{n}+1}^{n}\xi
_{ij}^{(n)}-\sum_{j=k_{n}+1}^{n}\zeta_{ij}^{(n)}\Biggr|
\geq\frac{\varepsilon
}{2}\sqrt{n\log n}\Biggr).\nonumber
\end{eqnarray}
Use the fact that $\sum_{j=k_{n}+1}^{n}\zeta_{ij}^{(n)}\sim\sqrt
{n-k_n}\cdot N(0, 1)$ and (\ref{zhengt}) to have
\begin{eqnarray*}
&&
P\Biggl(\sum_{j=k_{n}+1}^{n}\zeta_{ij}^{(n)}> \biggl(\beta-\frac{\varepsilon
}{2}\biggr)\sqrt{n\log n}\Biggr) \\
&&\qquad=  P \Biggl(N(0,1)> \biggl(\beta-\frac{\varepsilon
}{2}\biggr)\sqrt{\frac{n}{n-k_n}}\cdot\sqrt{\log n} \Biggr)\\
&&\qquad \geq\frac{C}{n^{(\beta-{\varepsilon}/{3})^2/2}\log n}
\end{eqnarray*}
uniformly for all $1\leq i \leq k_n$ as $n$ is sufficiently large and
as $0<\varepsilon/3<\beta$, where in the last inequality we use the fact
that $(\beta- (\varepsilon/2))\sqrt{n/(n-k_n)} \leq(\beta- (\varepsilon
/3))$ as $n$ is sufficiently large. This, (\ref{mue}) and (\ref
{bede}) imply
\begin{eqnarray*}
\max_{1\leq i \leq k_n}P\Biggl(\sum_{j=k_{n}+1}^{n}\xi_{ij}^{(n)}\leq
(\beta-\varepsilon)\sqrt{n\log n}\Biggr) &\leq& 1-\frac{C_1}{n^{(\beta
-{\varepsilon}/{3})^2/2}\log n} + \frac{C_2}{n^2(\log n)^3}\\
& \leq& 1-\frac{C_3}{n^{(\beta-{\varepsilon}/{3})^2/2}\log n}
\end{eqnarray*}
as $n$ is sufficiently large for any $0<\varepsilon/3 < \beta\leq2$.
Use inequality $1-x\leq e^{-x}$ for any $x>0$ to obtain
%
%
\begin{equation}\qquad
\max_{1\leq i\leq
k_{n}}P\Biggl(\sum_{j=k_{n}+1}^{n}\xi_{ij}^{(n)}\leq(\beta-\varepsilon
)\sqrt{n\log n}\Biggr)^{k_n} \leq\exp\bigl\{-Cn^{1-{(\beta
-\varepsilon/4)^2}/{2}} \bigr\}\hskip-12pt
\end{equation}
as $n$ is sufficiently large for any $0<\varepsilon/4 < \beta\leq2$.
From (\ref{meishu}), we conclude that
%
%
\begin{equation}\label{market}\quad
P\bigl(T_n\leq(\beta-2\varepsilon)\sqrt{n\log n}\bigr)\leq\exp\bigl\{
-Cn^{1-{(\beta-\varepsilon/4)^2}/{2}} \bigr\} + \frac{1}{n(\log n)^{3/2}}\hskip-12pt
\end{equation}
as $n$ is sufficiently large for any $0<\varepsilon/4 < \beta\leq2$.
Now, take $\beta=\sqrt{2}$, and we get
%
%
\begin{equation}\label{xiuqiu}
P \biggl(\frac{T_n}{\sqrt{n\log n}}\leq\sqrt{2}-2\varepsilon\biggr)
=O \biggl(\frac{1}{n(\log n)^{3/2}} \biggr)
\end{equation}
as $n\to\infty$ for sufficiently small $\varepsilon>0$. This and (\ref
{dedao}) imply (\ref{probasha}).\vspace*{10pt}

\textsc{Proof of (\ref{mada}) and (\ref{africa})}.\quad To prove these, it
suffices to show
%
%
\begin{equation}\label{sizzling}
\qquad
\limsup_{n\to\infty}\frac{T_n}{\sqrt{n\log n}}\leq2\qquad\mbox{a.s.}
\quad\mbox{and}\quad \liminf_{n\to\infty}\frac{T_n}{\sqrt{n\log n}}\geq
\sqrt{2}\qquad\mbox{a.s.}
\end{equation}
and
\begin{equation}
\label{blak}
P \biggl(\frac{T_n}{\sqrt{n\log n}}\in[a,b) \mbox{ for
infinitely many } n\geq2 \biggr)=1
\end{equation}
for any $(a,b)\subset(\sqrt{2}, 2)$.

First, choosing $\alpha=2$ in (\ref{lanse}), we have that $P(T_n\geq
(2+2\varepsilon)\sqrt{n\log n})=O(n^{-1}(\log n)^{-3})$ as $n\to\infty
$ for any $\varepsilon\in(0, 1)$. Thus, $\sum_{n\geq2}P(T_n\geq
(2+2\varepsilon)\times\break\sqrt{n\log n})<\infty$. By the Borel--Cantelli lemma,
\[
\limsup_{n\to\infty}\frac{T_n}{\sqrt{n\log
n}}\leq2+2\varepsilon\qquad\mbox{a.s.}
\]
for any $\varepsilon\in(0,1)$. This gives the first inequality in (\ref
{sizzling}). By the same reasoning, the second inequality follows from
(\ref{xiuqiu}). To prove (\ref{blak}), since $\{T_n, n\geq2\}$ are
independent from assumption, by the second Borel--Cantelli lemma, it is
enough to show
%
%
\begin{equation}\label{winn}
\sum_{n\geq2}P \biggl(\frac{T_n}{\sqrt{n\log n}}\in[a,b)
\biggr)=\infty
\end{equation}
for any $(a,b)\subset(\sqrt{2}, 2)$. By (\ref{lanse}), we have that
%
%
\begin{equation}\label{payton}
P \biggl(\frac{T_n}{\sqrt{n\log n}}\geq b \biggr)\leq\frac
{C}{n^{(b-\varepsilon)^2/2-1}}
\end{equation}
as $n$ is sufficiently large and $\varepsilon>0$ is sufficiently small.
By (\ref{most}),
\[
\max_{1\leq i\leq
k_{n}}\sum_{j=k_{n}+1}^{n}\xi_{ij}^{(n)} \leq T_{n} + V_n
\]
for $n\geq2$. Thus, by independence and (\ref{shengchu}),
%
%
\begin{eqnarray}\label{xiaonu}\qquad
&& P\bigl(T_n\geq a\sqrt{n\log n}\bigr)\nonumber\\
&&\qquad\geq P\Biggl(\max_{1\leq i\leq
k_{n}}\sum_{j=k_{n}+1}^{n}\xi_{ij}^{(n)}\geq(a+\varepsilon)\sqrt
{n\log n}\Biggr)\nonumber\\
&&\qquad\quad{} - P\bigl(V_n\geq\varepsilon\sqrt{n\log
n}\bigr)\\
&&\qquad\geq 1- \Biggl(1-\min_{1\leq i\leq
k_{n}}P\Biggl(\sum_{j=k_{n}+1}^{n}\xi_{ij}^{(n)}\geq(a+\varepsilon)\sqrt
{n\log n} \Biggr) \Biggr)^{k_n}\nonumber\\
&&\qquad\quad{} - \frac{1}{n(\log n)^{3/2}}\nonumber
\end{eqnarray}
as $n$ is sufficiently large. By (\ref{mue})
\begin{eqnarray*}
&& P\Biggl(\sum_{j=k_{n}+1}^{n}\xi_{ij}^{(n)}\geq(a+\varepsilon)\sqrt
{n\log n}\Biggr)\\
&&\qquad \geq P\Biggl(\sum_{j=k_{n}+1}^{n}\zeta_{ij}^{(n)}\geq(a+ 2\varepsilon
)\sqrt{n\log n}\Biggr)\\
&&\qquad\quad{} - P\Biggl(\Biggl|\sum_{j=k_{n}+1}^{n}\xi_{ij}^{(n)}-\sum
_{j=k_{n}+1}^{n}\zeta_{ij}^{(n)}\Biggr|\geq\varepsilon\sqrt{n\log n}\Biggr)\\
&&\qquad \geq P\bigl(N(0,1)\geq(a+ 3\varepsilon)\sqrt{\log n}\bigr)-\frac{1}{n^2}
\end{eqnarray*}
uniformly for all $1\leq i \leq k_n$ as $n$ is sufficiently large. From
(\ref{zhengt}), for any $\varepsilon>0$,
\[
P\bigl(N(0,1)\geq(a+ 3\varepsilon)\sqrt{\log n}\bigr)\sim\frac{C}{n^{(a+
3\varepsilon)^2/2}\sqrt{\log n}}
\]
as $n$ is sufficiently large. Noting that $a\in(\sqrt{2}, 2)$, we have
\[
P\Biggl(\sum_{j=k_{n}+1}^{n}\xi_{ij}^{(n)}\geq(a+\varepsilon)\sqrt{n\log
n}\Biggr)\geq\frac{C}{n^{(a+ 3\varepsilon)^2/2}\sqrt{\log n}}
\]
uniformly for all $1\leq i \leq k_n$ as $n$ is sufficiently large and
$\varepsilon$ is sufficiently small. Thus, since $k_n=[n/\log n]$, relate
the above to (\ref{xiaonu}) to give us that
\begin{eqnarray*}
P\bigl(T_n\geq a\sqrt{n\log n}\bigr) &\geq& 1- \biggl(1-\frac{C}{n^{(a+
3\varepsilon)^2/2}\sqrt{\log n}} \biggr)^{k_n}-\frac{1}{n(\log
n)^{3/2}}\\
& \sim& \frac{Ck_n}{n^{(a+ 3\varepsilon)^2/2}\sqrt{\log
n}}\bigl(1+o(1)\bigr)-\frac{1}{n(\log n)^{3/2}}\\
&\geq&\frac{C}{n^{(a+ 3\varepsilon
)^2/2-1}(\log n)^2}
\end{eqnarray*}
as $n$ is sufficiently large and $\varepsilon>0$ is small enough, where
in the ``$\sim$'' step above we use the fact that $1-(1-x_n)^{k_n}\sim
k_nx_n$ if $x_n\to0$, $k_n\to+\infty$ and $k_nx_n\to0$ as $n\to
\infty$. Combining this and (\ref{payton}), we eventually arrive at
\begin{eqnarray*}
P \biggl(\frac{T_n}{\sqrt{n\log n}}\in[a,b) \biggr)
&=& P\biggl(\frac{T_n}{\sqrt{n\log n}}\geq a \biggr)- P \biggl(\frac{T_n}{\sqrt
{n\log n}}\geq b \biggr)\\
&\geq& \frac{C_3}{n^{(a+ 3\varepsilon)^2/2-1}(\log n)^2}-\frac
{C_4}{n^{(b-\varepsilon)^2/2-1}}\\
&\sim& \frac{C_3}{n^{(a+ 3\varepsilon
)^2/2-1}(\log n)^2}
\end{eqnarray*}
as $n$ is sufficiently large and $\varepsilon>0$ is sufficiently small,
where $[a, b)\subset(\sqrt{2}, 2)$. Finally, choosing $\varepsilon>0$
so small that $(a+ 3\varepsilon)^2/2-1\in(0, 1)$, we get (\ref{winn}).
\end{pf*}
\begin{pf*}{Proof of Corollary \ref{coro1}} Recalling (\ref
{Mar}), let $\tilde{\xi}_{ij}^{(n)}=(\xi_{ij}^{(n)}-\mu_n)/\sigma
_n$ for all $1\leq i< j\leq n$ and $n\geq2$. Then $\{\tilde{\xi
}_{ij}^{(n)}; 1\leq i < j \leq n, n\geq2\}$ satisfies (\ref
{shufic}) with $
\mu_n=0, \sigma_n=1$ and $p>6$. Let $\tilde{\bd{\Delta}}_n$ be
generated by $\{\tilde{\xi}_{ij}^{(n)}\}$ as in (\ref{Mar}). By
Theorem~\ref{west1}, the conclusions there hold if $\lambda_{\max}(\bd
{\Delta}_n)$ is replaced by $\lambda_{\max}(\tilde{\bd{\Delta
}}_n)$. Notice
%
%
\begin{equation}\label{future}
\bd{\Delta}_n=\sigma_n\tilde{\bd{\Delta}}_n + \mu_n\cdot(n
\mathbf{I}_n - \mathbf{J}_n),
\end{equation}
where $\mathbf{I}_n$ is the $n\times n$ identity matrix, and $\mathbf{J}_n$
is the $n\times n$ matrix with all of its entries equal to 1. It is
easy to check that the eigenvalues of $n\mathbf{I}_n - \mathbf{J}_n$ are $0$
with one fold and $n$ with $n-1$ folds, respectively.\vspace*{1pt} First, apply the
triangle inequality to (\ref{future}) to have that $|\lambda
_{\max}(\bd{\Delta}_n)-\sigma_n\lambda_{\max}(\tilde{\bd{\Delta
}}_n)|\leq\|\mu_n\cdot(n\mathbf{I}_n - \mathbf{J}_n)\| \leq n|\mu_n|$. It
follows that
\[
\biggl|\frac{\lambda_{\max}(\bd{\Delta}_n)}{\sqrt{n\log n} \sigma
_n}-\frac{\lambda_{\max}(\tilde{\bd{\Delta}}_n)}{\sqrt{n\log
n}} \biggr|
\leq\frac{|\mu_n|}{(\log n)^{1/2}n^{-1/2}\sigma_n} \to0
\]
provided $|\mu_n| \ll\sigma_n\sqrt{\log n/n}$. Then (a1) and (b1)
follow from Theorem \ref{west1}. By the same argument
\[
\bigl|\lambda_{\max}(\bd{\Delta}_n) - \lambda_{\max}\bigl(\mu_n\cdot(n\mathbf{I}_n - \mathbf{J}_n)\bigr)\bigr|
\leq\sigma_n\|\tilde{\bd{\Delta}}_n\|=O\bigl(\sigma_n\sqrt{n\log
n}\bigr)\qquad
\mbox{a.s.}
\]
as $n\to\infty$. Note that $\lambda_{\max}(\mu_n\cdot(n\mathbf{I}_n -
\mathbf{J}_n))=0$ if $\mu_n<0$, and is equal to $n\mu_n$ if $\mu_n>0$
for any $n\geq2$. Thus, if $\mu_n\gg\sigma_n\sqrt{\log n/n}$, we
have $\lambda_{\max}(\bd{\Delta}_n)/(n\mu_n)\to1$ a.s. as $n\to
\infty$. If $\mu_n<0$ for all $n\geq2$, and $|\mu_n| \gg\sigma
_n\sqrt{\log n/n}$, we obtain $\lambda_{\max}(\bd{\Delta}_n)/(n\mu
_n)\to0$ a.s. as $n\to\infty$. Then (a2), (a3), (b2) and (b3) are yielded.

Finally, since $E(-\xi_{ij}^{(n)})=-\mu_n$ and $\operatorname{Var}(-\xi
_{ij}^{(n)})=\operatorname{Var}(\xi_{ij}^{(n)})=\sigma_n^2$ for all $i,j,n$, by
using the proved (a1) and (b1), we know that (a1) and (b1) are also
true if $\lambda_{\max}(\bd{\Delta}_n)$ is replaced by $\lambda
_{\max}(-\bd{\Delta}_n)$. Now, use the same arguments as in the proof
of part (c) in Theorem \ref{west1} to get (a1) and (b1) when $\lambda
_{\max}(\bd{\Delta}_n)$ is replaced with $\|\bd{\Delta}_n\|$. On the
other hand, it is well known that $\bd{\Delta}_n$ is nonnegative
definite if $\xi_{ij}^{(n)}\geq0$ for all $i,j,n$ (see, e.g., page 5 in
\cite{Chung97}). Thus $\|\bd{\Delta}_n\|=\lambda_{\max}(\bd{\Delta
}_n)$. Consequently (a2) and (b2) follow when $\lambda_{\max}(\bd
{\Delta}_n)$ is replaced with $\|\bd{\Delta}_n\|$.
\end{pf*}

To prove Theorem \ref{free1}, we need some preliminary results.
\begin{lemma}\label{whoknows}
Let $\{\xi_{ij}^{(n)}; 1\leq i<j\leq n, n\geq2\}$ be defined on the
same probability space. For each $n\geq2$, let $\{\xi_{ij}^{(n)};
1\leq i<j\leq n\}$ be
independent r.v.s with $E\xi_{ij}^{(n)}=0$. Define $\xi
_{ji}^{(n)}=\xi_{ij}^{(n)}$ for all $i,j,n$ and $S_{n,1}=\sum_{1\leq
i\ne j \leq n}(\xi_{ij}^{(n)})^2$ and $S_{n,2}=\sum_{i=1}^{n}(\sum
_{j\neq
i}\xi_{ij}^{(n)})^{2}$. If $\sup_{1\leq
i<j\leq n, n\geq2}E|\xi_{ij}^{(n)}|^{4+\delta}<\infty$ for some
$\delta>0$, then
%
%
\begin{equation}\label{nou}
\lim_{n\to\infty}\frac{S_{n, k} - ES_{n,k}}{n^{2}}=0 \qquad\mbox{a.s.}
\qquad\mbox{for } k=1,2.
\end{equation}
\end{lemma}
\begin{pf}
To make notation simple, we write $\xi_{ij}=\xi
_{ij}^{(n)}$ for all $1\leq i\leq j\leq n$ when there is no confusion.

\textit{Case} 1: $k=1$. Recall the Marcinkiewicz--Zygmund inequality (see,
e.g., Corollary 2 and its proof on page 368 in \cite{ChowT}), for any
$p\geq2$, there exists a constant $C_p$ depending on $p$ only such that
%
%
\begin{equation}\label{MZ}
E\Biggl|\sum_{i=1}^nX_i\Biggr|^p\leq C_p n^{{p}/{2}-1} \sum_{i=1}^nE|X_i|^p
\end{equation}
for any sequence of independent random variables $\{X_i; 1\leq i \leq
n\}$ with $EX_i=0$ and $E(|X_i|^p)< \infty$ for all $1\leq i\leq n$.
Taking $p=2+(\delta/2)$ in (\ref{MZ}), we have from the H\"{o}lder
inequality that
%
%
\begin{eqnarray}\label{wen1}
E(|\xi_{ij}^2-E\xi_{ij}^2|^{p}) & \leq& 2^{p-1} E|\xi_{ij}|^{2p} +
2^{p-1}(E|\xi_{ij}|^{2})^{p}\nonumber\\
&\leq& 2^p\cdot\sup_{1\leq
i,j\leq n,
n\geq1}E\bigl|\xi_{ij}^{(n)}\bigr|^{4+\delta}\\
&<&\infty\nonumber
\end{eqnarray}
uniformly for all $1\leq
i<j\leq n, n\geq2$. Write $S_{n,1}-ES_{n,1}=2\sum_{1\leq i < j\leq
n}(\xi_{ij}^2-E\xi_{ij}^2)$. By (\ref{MZ}),
%
%
\begin{eqnarray}\label{lou2}\quad
E|S_{n,1}-ES_{n,1}|^{p}&\leq& C\cdot\biggl(\frac{n(n-1)}{2}
\biggr)^{\delta/4}\cdot\sum_{1\leq i < j\leq n}
E(|\xi_{ij}^2-E\xi_{ij}^2|^{p})\nonumber\\[-8pt]\\[-8pt]
& \leq& C\cdot n^{2+(\delta/2)},\nonumber
\end{eqnarray}
where $C$ here and later, as earlier, is a constant not depending on
$n$, and may be different from line to line. Then
$P(|S_{n,1}-ES_{n,1}|\geq n^2\varepsilon)\leq(n^2\varepsilon
)^{-p}E|S_{n,1}-ES_{n,1}|^p=O(n^{-2-(\delta/2)})$ for any $\varepsilon
>0$ by the Markov inequality. Then (\ref{nou}) holds for $k=1$ by the
Borel--Cantelli lemma.

\textit{Case} 2: $k=2$. For $n\geq2$, set $u_1=v_n=0$ and
\[
u_i=\sum_{j=1}^{i-1}\xi_{ij} \qquad\mbox{for } 2\leq i \leq n+1
\quad\mbox{and}\quad
v_i=\sum_{j=i+1}^n\xi_{ij} \qquad\mbox{for } 0\leq i \leq n-1.
\]
Then, $\sum_{j\neq i}\xi_{ij}=u_i+v_i$ for all $1\leq i \leq n$.
%
%
%
Clearly, $S_{n,2}=\sum_{i=1}^nu_i^2+ \sum_{i=1}^nv_i^2+ 2\sum
_{i=1}^nu_iv_i$ for all $n\geq1$. Since $E(u_iv_i)=(Eu_i) Ev_i=0$ by
independence, to prove the lemma for $k=2$, it suffices to show
\[
\frac{1}{n^2}\sum_{i=1}^n(u_i^2-Eu_i^2)\to0 \qquad\mbox{a.s.},\qquad
\frac{1}{n^2}\sum_{i=1}^n(v_i^2-Ev_i^2)\to0\qquad\mbox{a.s.}
\]
%
%
and
\begin{equation}\label{smoke}
\frac{1}{n^2}\sum_{i=1}^nu_iv_i\to0\qquad\mbox{a.s.}
\end{equation}
as $n\to\infty$. We will only prove the first and the last assertions
in two steps. The proof of the middle one is almost the same as that of
the first and, therefore, is omitted.

\textit{Step} 1. 
Similarly\vspace*{1pt} to the discussion in (\ref{wen1}) and (\ref{lou2}), we have
$E|u_{i}|^{4+\delta}\leq C i^{2+(\delta/2)}$ for all $1\leq i \leq
n$ and $n\geq2$. Now set $Y_{n,i}=(u_{i}^2-Eu_i^2)/i$ for
$i=1,2,\ldots, n$. Then, $\{Y_{n,i}; 1\leq i \leq n\}$ are
independent random variables with
%
%
\begin{equation}
\label{nde}
EY_{n,i}=0,\qquad \sup_{1\leq
i,j\leq n, n\geq1}E|Y_{n,i}|^{2+\delta'}< \infty
\end{equation}
and
\begin{equation}
\label{wode}
\frac{1}{n^2}\sum_{i=1}^n(u_i^2-Eu_i^2)=\frac{1}{n^2}\sum
_{i=1}^{n}iY_{n,i}
\end{equation}
for all $1\leq i \leq n$ and $n\geq2$, where $\delta'=\delta/2$. By
(\ref{MZ}) and (\ref{nde}),
\[
E\Biggl|\sum_{i=1}^{n}iY_{n,i}\Biggr|^{2+\delta'}\leq C\cdot n^{(2+\delta
')/2-1}\sum_{i=1}^{n}i^{2+\delta'}=O\bigl(n^{3+(3\delta'/2)}\bigr)
\]
as $n\to\infty$, where the inequality $\sum_{i=1}^{n}i^{2+\delta
'}\leq\sum_{i=1}^{n}n^{2+\delta'}\leq n^{3+\delta'}$ is used in the
above inequality. For any $t>0$,
\[
P \Biggl(\frac{1}{n^2}\Biggl|\sum_{i=1}^{n}iY_{n,i}\Biggr|\geq t \Biggr)\leq\frac
{E|\sum_{i=1}^{n}iY_{n,i}|^{2+\delta'}}{(n^2t)^{2+\delta'}}
=O \biggl(\frac{1}{n^{1+(\delta'/2)}} \biggr)
\]
as $n\to\infty$. This together with (\ref{wode}) concludes the first
limit in (\ref{smoke}) by the Borel--Cantelli lemma.

\textit{Step} 2. We will prove the last assertion in (\ref{smoke}) in
this step.
Define $\sigma$-algebra
\[
\mathcal{F}_{n,0}=\{\varnothing, \Omega\} \quad\mbox{and}\quad \mathcal
{F}_{n,k}=\sigma\bigl(\xi_{ij}^{(n)}; 1\leq i\leq k, i+1\leq j
\leq n \bigr)
\]
for $1\leq k \leq n-1$. Obviously, $\mathcal{F}_{n,0}\subset\mathcal
{F}_{n,1}\subset\cdots\subset\mathcal{F}_{n,n-1}$. It is easy to
verify that
\[
E \Biggl( \sum_{i=1}^{k+1}u_iv_i | \mathcal{F}_{n,k} \Biggr)=\sum
_{i=1}^{k}u_iv_i
\]
for $k=1,2,\ldots, n-1$. Therefore, $\{\sum_{i=1}^{k}u_iv_i,
\mathcal{F}_{n,k}, 1\leq k\leq n-1\}$ is a martingale. By the given
moment condition, $\tau:=\sup_{1\leq
i,j\leq n,n\geq1}E|\xi_{ij}^{(n)}|^{4}<\infty$. From (\ref{MZ}),
$E(u_i^4)\leq Ci^{2}\leq Cn^2$ and $E(v_i^4)\leq C(n-i)^{2}\leq Cn^2$
for $1\leq i \leq n$ and $n\geq2$. By applying the Burkholder
inequality (see, e.g., Theorem 2.10 from
\cite{HH} or Theorem~1 on page 396 and the proof of Corollary 2 on
page 268 from \cite{ChowT}), we have
\[
E \Biggl(\sum_{i=1}^{n-1}u_iv_i \Biggr)^4 \leq Cn^{(4/2)-1}\sum
_{i=1}^{n-1}E((u_iv_i)^4)=Cn\sum_{i=1}^{n-1}E(u_i)^4\cdot E(v_i)^4=O(n^6)
\]
as $n\to\infty$.
By the Markov inequality,
\[
P \Biggl(\frac{1}{n^2}\Biggl|\sum_{i=1}^nu_iv_i\Biggr|\geq\delta\Biggr)\leq
\frac{E|\sum_{i=1}^{n-1}u_iv_i|^4}{n^8\delta^4}=O \biggl(\frac
{1}{n^2} \biggr)
\]
as $n\to\infty$. The Borel--Cantelli says that $\sum
_{i=1}^nu_iv_i/n^2\to0$ a.s. as $n\to\infty$.
\end{pf}

For any two probability measures $\mu$ and $\nu$ on $\mathbb{R}$, define
%
%
\begin{equation}
d_{\mathrm{BL}}(\mu, \nu)=\sup\biggl\{\int f \,d\mu-\int f \,d\nu
\dvtx\|f\|_{\infty}+\|f\|_{L}\leq1\biggr\},
\end{equation}
where $\|f\|_{\infty}=\sup_{x\in\mathbb{R}}|f(x)|, \|f\|_{L}=\sup
_{x\neq
y}|f(x)-f(y)| /|x-y|$. It is well known (see, e.g., Section 11.3 from
\cite{Dudley}),
that $d_{\mathrm{BL}}(\cdot, \cdot)$ is called the bounded Lipschitz metric,
which characterizes the weak convergence of probability measures.
Reviewing (\ref{minyi}), for
the spectral measures of $n\times n$ real and symmetric matrices
$\mathbf{M}_1$ and $\mathbf{M}_2$, we
have (see, e.g., (2.16) from \cite{BDJ})
%
%
\begin{equation}\label{judy}
d_{\mathrm{BL}}^{2}(\hat{\mu}(\mathbf{M}_1),\hat{\mu}(\mathbf{M}_2))\leq\frac
{1}{n}\operatorname{tr}\bigl((\mathbf{M}_1-\mathbf{M}_2)^{2}\bigr).
\end{equation}

To prove Theorem \ref{free1}, we first reduce it to the case that all
random variables in the matrices are uniformly
bounded. This step will be carried out through a truncation argument by
using (\ref{judy}).
\begin{lemma}\label{sunnyday}
If Theorem \ref{free1} holds for all uniformly bounded
r.v.s $\{\xi_{ij}^{(n)};1\leq i<j\leq n, n\geq2\}$ satisfying (\ref
{shufic}) with $\mu_n=0$ and $\sigma_n=1$ for all $n\geq2$, then it
also holds for all
r.v.s $\{\xi_{ij}^{(n)};1\leq i<j\leq n, n\geq2\}$ satisfying (\ref
{shufic}) with $p=4+\delta$ for some $\delta>0$, and $\mu_n=0$ and
$\sigma_n=1$ for all $n\geq2$.
\end{lemma}
\begin{pf}
As in the proof of Lemma \ref{whoknows}, we write
$\xi_{ij}$ for $\xi_{ij}^{(n)}$ if there is no danger of confusion.
Fix $u>0$. Let
\[
\tilde{\xi}_{ij}=\xi_{ij}I\{|\xi_{ij}|\leq u\} - E(\xi_{ij}I\{
|\xi_{ij}|\leq u\})
\]
and
\[
\sigma_{ij}(u)=\sqrt{\operatorname{Var}(\tilde{\xi}_{ij})}
\]
for all $i$ and $j$. Note that
\[
\bigl|\sigma_{ij}(u)-\sqrt{\operatorname{Var}(\xi
_{ij})}\bigr|\leq\sqrt{\operatorname{Var}(\xi_{ij}-\tilde{\xi}_{ij})}\leq\sqrt{E\xi
_{ij}^2I\{|\xi_{ij}|> u\}}
\]
by the triangle inequality. Thus, with
condition that $\sup_{1\leq i<j\leq n, n\geq2}E|\xi
_{ij}^{(n)}|^{4+\delta}<\infty$, we see that
%
%
\begin{equation}\label{howk}\qquad
{\sup_{1\leq i<j\leq n, n\geq2}}|\sigma_{ij}(u)-1|\to0
\quad\mbox{and}\quad
\sup_{1\leq i<j\leq n, n\geq2}E(\xi_{ij}-\tilde{\xi}_{ij})^2\to0
\end{equation}
as $u\to+\infty$. Take $u>0$ large enough such that $\sigma
_{ij}(u)>1/2$ for all $1\leq i\ne j \leq n$ and $n\geq2$. Write
\[
\xi_{ij}=\underbrace{\frac{\tilde{\xi}_{ij}}{\sigma
_{ij}(u)}}_{x_{ij}^{(n)}}{} + {}\underbrace{ \frac{\sigma
_{ij}(u)-1}{\sigma_{ij}(u)}\cdot\tilde{\xi}_{ij}}_{y_{ij}^{(n)}}{}
+{}\underbrace{(\xi_{ij}-\tilde{\xi}_{ij})}_{z_{ij}^{(n)}}
\]
for all $1\leq i\neq j \leq n, n\geq2$. Obviously, for
$a_{ij}^{(n)}=x_{ij}^{(n)}, y_{ij}^{(n)}$ or $z_{ij}^{(n)}$, we know
$\{a_{ij}^{(n)}; 1\leq i< j\leq n\}$ are independent for each $n\geq
2$, and
%
%
\begin{equation}\label{breada}
Ea_{ij}^{(n)}=0 \quad\mbox{and}\quad \sup_{1\leq i<j\leq n, n\geq
2}E\bigl|a_{ij}^{(n)}\bigr|^{4+\delta}<\infty.
\end{equation}
Again, for convenience, write $x_{ij}, y_{ij}$ and $z_{ij}$ for
$x_{ij}^{(n)}, y_{ij}^{(n)}$ and $z_{ij}^{(n)}$. Clearly, $\{
x_{ij}; 1\leq i < j\leq n, n\geq2\}$ are uniformly bounded.
Besides, it is easy to see from (\ref{howk}) that
%
%
\begin{equation}\label{aawa}
\sup_{1\leq i<j\leq n, n\geq2}\bigl(E(y_{ij}^2) + E(z_{ij}^2)\bigr) \to0
\end{equation}
as $u\to+\infty$.

Let $\mathbf{X}_n, \mathbf{Y}_n$ and $\mathbf{Z}_n$ be the Laplacian matrices
generated by $\{x_{ij}\}, \{y_{ij}\}$ and $\{z_{ij}\}$ as in (\ref
{Mar}), respectively. Then $\bd{\Delta}_n=\mathbf{X}_n +\mathbf{Y}_n +
\mathbf{Z}_n$. With (\ref{judy}), use the inequality that $\operatorname{tr}((\mathbf{M}_1+
\mathbf{M}_2)^2)\leq2 \operatorname{tr}(\mathbf{M}_1^2) +2 \operatorname{tr}(\mathbf{M}_2^2)$ for any
symmetric matrices $\mathbf{M}_1$ and $\mathbf{M}_2$ to obtain that
\begin{eqnarray*}
d_{\mathrm{BL}}^2 \biggl(\frac{\bd{\Delta}_n}{\sqrt{n}},
\frac{\mathbf{X}_{n}}{\sqrt{n}} \biggr)
& \leq& \frac{1}{n^2}\operatorname{tr}\bigl((\mathbf{Y}_n + \mathbf{Z}_n)^2\bigr)\\
&\leq& \frac{2}{n^2}\sum_{1\leq i\ne j\leq n}\bigl((y_{ij})^2+
(z_{ij})^2\bigr)\\
&&{} +\frac{2}{n^2}\sum_{i=1}^n\biggl\{\biggl(\sum_{j\ne i}y_{ij}\biggr)^2 +
\biggl(\sum_{j\ne i}z_{ij}\biggr)^2\biggr\}.
\end{eqnarray*}
By independence and symmetry,
\[
E\biggl(\biggl(\sum_{j\ne i}y_{ij}\biggr)^2 + \biggl(\sum
_{j\ne i}z_{ij}\biggr)^2\biggr)=2\sum_{j\ne i}\{E(y_{ij})^2+E(z_{ij})^2\}.
\]
Recalling (\ref{breada}), by applying Lemma \ref{whoknows}, we have
%
%
\begin{eqnarray}\label{eveny}
&&\limsup_{n\to\infty}d_{\mathrm{BL}}^2 \biggl(\frac{\bd{\Delta}_n}{\sqrt
{n}}, \frac{\mathbf{X}_{n}}{\sqrt{n}} \biggr)\nonumber\\[-8pt]\\[-8pt]
&&\qquad\leq C\cdot\sup_{1\leq
i<j\leq n, n\geq2}\bigl(E(y_{ij}^2) +
E(z_{ij}^2)\bigr)\to0\qquad\mbox{a.s.}\nonumber
\end{eqnarray}
as $u\to+\infty$ thanks to (\ref{aawa}). Noticing $Ex_{ij}=0,
Ex_{ij}^2=1$ for all $i,j$, and $\{x_{ij}; 1\leq i< j \leq n, n\geq
2\}$ are uniformly bounded. By assumption,\break $d_{\mathrm{BL}}(\hat{\mu
}(n^{-1/2} 
\mathbf{X}_{n}), \gamma_{M})\to0$ as $n\to\infty$, where
$\gamma_{M}$ is the probability measure mentioned in Theorem \ref
{free1}. With this, (\ref{eveny}) and the triangle inequality of
metric $d_{\mathrm{BL}}$, we see that $d_{\mathrm{BL}}(\hat{\mu}(n^{-1/2}\bd{\Delta
}_{n}), \gamma_{M})\to0$ as $n\to\infty$.
\end{pf}



Given $n\geq2$, let $\Gamma_n=\{(i, j); 1\leq j < i\leq n\}$ be a
graph. We say $a=(i_1, j_1)$ and $b=(i_2, j_2)$ form an edge and denote
it by $a \sim b$, if one of $i_1$ and $j_1$ is identical to one of
$i_2$ and $j_2$. For convenience of notation, from now on, we write
$a=(a^+, a^-)$ for any $a\in\Gamma_n$. Of course, $a^+>a^-$. Given
$a, b\in\Gamma_n$, define an $n\times n$ matrix
\[
\mathbf{Q}_{a,b}[i,j]
= \cases{
-1,&\quad if $ i=a^{+}, j=b^{+}$ or $ i=a^{-}, j=b^{-}$;\cr
1, &\quad if $ i=a^{+}, j=b^{-}$ or $ i=a^{-}, j=b^{+}$;\cr
0, &\quad otherwise.}
\]
With this notation, we rewrite $\mathbf{M}_n$ as follows
%
%
\begin{equation}\label{manh}
-\bd{\Delta}_{n}=\sum_{a\in\Gamma_{n}}\xi_{a}^{(n)}\mathbf{Q}_{a,a},
\end{equation}
where $\xi_{a}^{(n)}=\xi_{a^+a^-}^{(n)}$ for $a\in\Gamma_n$. Let
$t_{a,b}=\operatorname{tr}(\mathbf{Q}_{a,b})$. We summarize some facts from \cite{BDJ}
in the following lemma.
\begin{lemma}\label{BDJ} Let $a, b\in\Gamma_n$. The following
assertions hold:

\begin{longlist}
\item $t_{a, b}=t_{b, a}$.\vspace*{28pt}

\item
\begin{eqnarray*}
\\[-58pt]
t_{a,b} &=&\cases{
-2, &\quad if $ a=b$;\cr
-1, &\quad if $ a\neq b$ and $ a^{-}=b^{-}$ or $a^{+}=b^{+};$\cr
1, &\quad if $ a\neq b$ and $a^{-}=b^{+}$ or $a^{+}=b^{-}$;\cr
0, &\quad otherwise.}
\end{eqnarray*}

\item $\mathbf{Q}_{a,b}\times
\mathbf{Q}_{c,d}=t_{b,c}\mathbf{Q}_{a,d}$. Therefore, $\operatorname{tr}(
\mathbf{Q}_{a_{1},a_{1}}\times\mathbf{Q}_{a_{2},a_{2}}\times\cdots\times
\mathbf{Q}_{a_{r},a_{r}})=\prod_{j=1}^{r}t_{a_{j},a_{j+1}}$, where $a_1,
\ldots, a_r\in\Gamma_n$, and $a_{r+1}=a_{1}$.
\end{longlist}
\end{lemma}

We call $\pi=(a_{1}, \ldots, a_{r})$ a circuit of length $r$ if
$a_{1}\sim\cdots\sim a_{r}\sim a_{1}$. For such a circuit, let
%
%
\begin{equation}\label{froa}
\xi_{\pi}^{(n)}=\prod_{j=1}^{r}t_{a_{j},a_{j+1}}\prod_{j=1}^{r}\xi
_{a_{j}}^{(n)}.
\end{equation}
From (\ref{manh}), we know
%
%
\begin{equation}\label{dingz}
\operatorname{tr}(\bd{\Delta}_{n}^{r})=(-1)^r\sum_{\pi}\xi_{\pi}^{(n)} \quad
\mbox{and}\quad E
\operatorname{tr}(\bd{\Delta}_{n}^{r})=(-1)^r\sum_{\pi}E\xi_{\pi}^{(n)},
\end{equation}
where the sum is taken over all circuits of length $r$ in $\Gamma_{n}$.
\begin{definition} We say that a circuit $\pi=(a_{1}\sim\cdots\sim
a_{r}\sim a_{1})$ of
length $r$ in $\Gamma_{n}$ is vertex-matched if for each
$i=1,\ldots, r$ there exists some $j\neq i$ such that $a_{i}=a_{j}$,
and that it has a match of order $3$ if some value is repeated at
least three times among $\{a_{j},j=1,\ldots, r\}$.
\end{definition}

Clearly, by independence, the only possible
nonzero terms in the second sum in (\ref{dingz}) come from
vertex-matched circuits. For $x\geq0$, denote by $\lfloor x\rfloor$
the integer part of $x$. The following two lemmas will be used later.
\begin{lemma}[(Propositions 4.10 and 4.14 from \cite{BDJ})]\label{BDJprop1}
Fix $r\in\mathbb{N}$.

\begin{longlist}
\item Let N denote the number of vertex-matched
circuits in $\Gamma_{n}$ with vertices having at least one match
of order 3. Then $N=O(n^{\lfloor(r+1)/2\rfloor})$ as $n\to\infty$.

\item Let N denote the number of vertex-matched
quadruples of circuits in $\Gamma_{n}$ with r vertices each, such
that none of them is self-matched. Then $N=O(n^{2r+2})$ as $n\to\infty$.
\end{longlist}
\end{lemma}

Let $\mathbf{U}_n$ be a symmetric matrix of form
%
%
\begin{equation}\label{heiy}
\mathbf{\mathbf{U}_n}= \pmatrix{
\displaystyle\sum_{j\neq1}Y_{1j}&-Y_{12}&-Y_{13}&\cdots&-Y_{1n}\vspace*{2pt}\cr
-Y_{21}&\displaystyle\sum_{j\neq2}Y_{2j}&-Y_{23}&\cdots&-Y_{2n}\vspace*{2pt}\cr
-Y_{31}&-Y_{32}&\displaystyle\sum_{j\neq3}Y_{3j}&\cdots&-Y_{3n}\vspace*{2pt}\cr
\vdots&\vdots&\vdots&\vdots&\vdots\vspace*{2pt}\cr
-Y_{n1}&-Y_{n2}&-Y_{n3}&\cdots&\displaystyle\sum_{j\neq n}Y_{nj}},
\end{equation}
where $\{Y_{ij};1\leq i< j<\infty\}$ are
i.i.d. standard normal random variables not depending on $n$.
\begin{lemma}\label{odd1} Suppose the conditions in Theorem \ref{free1}
hold with $\mu_n=0$ and $\sigma_n=1$ for all $n\geq2$. Furthermore,
assume $\{\xi_{ij}^{(n)}; 1\leq i <j \leq n, n\geq2\}$ are
uniformly bounded. Then:
\begin{longlist}
\item
$\lim_{n\to\infty}\frac{1}{n^{k+1/2}}\mathbb{E}\operatorname{tr}(\bd
{\Delta}_{n}^{2k-1})=0$;

\item
$\lim_{n\to\infty}\frac{1}{n^{k+1}} (\mathbb{E}\operatorname{tr}(\bd
{\Delta}_{n}^{2k})-\mathbb{E}\operatorname{tr}(\mathbf{U}_n^{2k}) )=0$
\end{longlist}
for any integer $k\geq1$, where $\mathbf{U}_n$ is as in (\ref{heiy}).
\end{lemma}
%
%
\begin{pf}
(i) As remarked earlier, all nonvanishing
terms in the representation of $\mathbb{E}\operatorname{tr}(\bd{\Delta
}_{n}^{2k-1})$ in (\ref{dingz}) are of form $\mathbb{E}\xi_{\pi
}^{(n)}$ with the vertices of the path $a_{1}\sim a_{2}\sim\cdots\sim
a_{2k-1}\sim a_1$ in $\pi$ repeating at least two times. Since $2k-1$
is an odd number, there exists a vertex such that it repeats at least
three times. Also, in view of (\ref{froa}) and that $|t_{a,b}|\leq2$
for any $a, b \in\Gamma_n$, thus all such terms $\mathbb{E}\xi_{\pi
}^{(n)}$ are uniformly bounded.
Therefore, by (i) of Lemma \ref{BDJprop1},
\[
\biggl|\frac{1}{n^{k+1/2}}\cdot\mathbb{E}\operatorname{tr}(\bd{\Delta}_{n}^{2k-1})\biggr|\leq
\frac{C}{\sqrt{n}}\to0
\]
as $n\to\infty$, where $C$ is a constant not depending on $n$.

(ii) Recall (\ref{heiy}). Define $Y_{\pi}^{(n)}$ similarly to $\xi
_{\pi}^{(n)}$ in (\ref{froa}). We then have that
\begin{eqnarray*}
|\mathbb{E}\operatorname{tr}(\bd{\Delta}_{n}^{2k})-\mathbb{E}\operatorname{tr}(
\mathbf{U}_{n}^{2k})|&=&
\biggl|\sum_{\pi}\bigl(\mathbb{E}\xi_{\pi}^{(n)}-\mathbb{E}Y_{\pi
}^{(n)}\bigr)\biggr|\\
&\leq&\biggl|\sum_{\pi\in
A_{1}}\bigl(\mathbb{E}\xi_{\pi}^{(n)}-\mathbb{E}Y_{\pi}^{(n)}\bigr) \biggr|+
\biggl|\sum_{\pi\in A_{2}}\bigl(\mathbb{E}\xi_{\pi}^{(n)}-\mathbb
{E}Y_{\pi}^{(n)}\bigr)
\biggr|\\
&:=&I_{1}+I_{2},
\end{eqnarray*}
where $A_{1}$ denotes the set of the vertex-matched circuits with match
of order $3$, and
$A_{2}$ denotes the set of the vertex-matched circuits in $\Gamma_{n}$
such that there are exactly $k$ distinct matches. Observe that each
vertex of any circuit in $A_2$ matches exactly two times. From the
independence assumption and that $E|\xi_{ij}^{(n)}|^2=1$ for all
$1\leq i <j\leq n$ and $n\geq2$, we know $\mathbb{E}\xi_{\pi
}^{(n)}=\mathbb{E}Y_{\pi}^{(n)}=1$ for $\pi\in A_2$. This gives
$I_2=0$. By Lemma \ref{BDJprop1}, the cardinality of $A_{1}\leq
n^{k}$. Since $\xi_{ij}^{(n)}$ are uniformly
bounded and $Y_{ij}$ are standard normal random variables, we have
$I_{1}\leq
C n^{k}$ for some constant $C>0$ not depending on $n$. In summary
\[
\frac{1}{n^{k+1}}|\mathbb{E}\operatorname{tr}(\bd{\Delta}_{n}^{2k})-\mathbb
{E}\operatorname{tr}(\mathbf{U}_n^{2k})|=O \biggl(\frac{1}{n} \biggr)
\]
as $n\to\infty$. The proof is complete.
\end{pf}
\begin{lemma}\label{hahate} Suppose (\ref{shufic}) holds for some
$p>4$. Assume $\mu_n=0, \sigma_n=1$ for all $n\geq2$. Then, as
$n\to\infty$, $F^{\bd{\Delta}_n/\sqrt{n}}$ converges weakly
to the free convolution $\gamma_{M}$ of the
semicircular law and standard normal distribution. The measure
$\gamma_{M}$
is a nonrandom symmetric probability measure with
smooth bounded density, does not depend on the distribution of
$\{\xi_{ij}^{(n)}; 1\leq i<j\leq n, n\geq2\}$ and has an unbounded
support.
\end{lemma}
\begin{pf}
By Lemma \ref{sunnyday}, without loss of generality,
we now assume that $\{\xi_{ij}^{(n)}; 1\leq i < j \leq n, n\geq2\}
$ are uniformly bounded random variables with mean zero and variance
one, and $\{\xi_{ij}^{(n)}; 1\leq i < j \leq n\}$ are independent
for each $n\geq2$.

Proposition A.3 from \cite{BDJ} says that $\gamma_{M}$ is a symmetric
distribution and uniquely determined by its
moments. Thus, to prove the theorem, it is enough to show that
%
%
\begin{eqnarray}\label{hainr}
\frac{1}{n}\operatorname{tr}(n^{-1/2}\bd{\Delta}_{n})^{r}&=&\frac
{1}{n^{r/2+1}}\operatorname{tr}(\bd{\Delta}_{n}^{r})=\int
x^{r}\,dF^{n^{-1/2}\bd{\Delta}_{n}} \nonumber\\[-8pt]\\[-8pt]
& \to&\int
x^{r}\,d\gamma_{M} \qquad\mbox{as } n\to\infty
\qquad\mbox{a.s.}\nonumber
\end{eqnarray}
for any integer $r\geq1$. First, we claim that
%
%
\begin{equation}\label{wusi}
\mathbb{E}\bigl[\bigl(\operatorname{tr}(\bd{\Delta}_{n}^{r})-\mathbb{E}\operatorname{tr}(\bd{\Delta
}_{n}^{r})\bigr)^{4}\bigr]= O(n^{2r+2})
\end{equation}
as $n\to\infty$.
In fact, by (\ref{dingz}), we have
%
%
\begin{equation}\label{dita}
\mathbb{E}\bigl[\bigl(\operatorname{tr}(\bd{\Delta}_{n}^{r})-\mathbb{E}\operatorname{tr}(\bd{\Delta
}_{n}^{r})\bigr)^{4}\bigr]=\sum_{\pi_{1},\pi_{2},\pi_{3},\pi_{4}}
\mathbb{E}\Biggl[\prod_{j=1}^{4}\bigl(\xi_{\pi_{j}}-\mathbb{E}(\xi_{\pi_{j}})\bigr)\Biggr],
\end{equation}
where the sum runs over all circuits $\pi_{j}, j=1, 2, 3 ,4$ in
$\Gamma_{n}$, each having $r$ vertices. From the assumption, we know
$\{\xi_{ij}^{(n)},1\leq i<j\leq n\}$ are independent random variables
of mean
zero, and it is enough to consider the terms in (\ref{dita}) with all
vertex-matched
quadruples of circuits on $\Gamma_{n}$, such that none of them is
self-matched. By assumption, $\{\xi_{ij}^{(n)}; 1\leq i< j\leq n;
n\geq2\}$ are uniformly bounded, so all terms $\mathbb{E}[\prod
_{j=1}^{4}(\xi_{\pi_{j}}-\mathbb{E}(\xi_{\pi_{j}}))]$ in the sum
are uniformly bounded. By (ii) of Lemma \ref{BDJprop1}, we obtain
(\ref{wusi}).

By the Markov inequality,
%
%
\begin{eqnarray}
&& P \biggl(\frac{1}{n}|\operatorname{tr}((n^{-1/2}\bd{\Delta}_{n})^{r})-\mathbb
{E}\operatorname{tr}((n^{-1/2}\bd{\Delta}_{n})^{r})|\geq\varepsilon\biggr)
\nonumber\\[-8pt]\\[-8pt]
&&\qquad \leq
\frac{E|\operatorname{tr}(\bd{\Delta}_{n}^{r})-\mathbb{E}\operatorname{tr}(\bd{\Delta
}_{n}^{r})|^4}{(n^{1+(r/2)}\varepsilon)^4}=O \biggl(\frac{1}{n^{2}}
\biggr)\nonumber
\end{eqnarray}
as $n\to\infty$. It follows from the Borel--Cantelli lemma that
%
%
\begin{equation}\label{shiyan}
\frac{1}{n} \bigl(\operatorname{tr}((n^{-1/2}\bd{\Delta}_{n})^{r})-\mathbb
{E}\operatorname{tr}((n^{-1/2}\bd{\Delta}_{n})^{r}) \bigr)\to0\qquad\mbox{a.s.}
\end{equation}
as $n\to\infty$. Recalling $\mathbf{U}_n$ in (\ref{heiy}), Proposition
4.13 in \cite{BDJ} says that
\[
\frac{1}{n}\mathbb{E} \operatorname{tr}((n^{-1/2}\mathbf{U}_{n})^{2k}) \to\int
_{\mathbb{R}}x^{2k}\,d \gamma_{M}
\]
as $n\to\infty$ for any $k\geq1$. This, (ii) of Lemma \ref{odd1}
and (\ref{shiyan}) imply (\ref{hainr}) for any even number $r\geq1$.
For odd number $r$, (i) of Lemma \ref{odd1} and (\ref{shiyan}) yield
(\ref{hainr}) since $\gamma_{M}$ is symmetric, hence its odd moments
are equal to zero.
\end{pf}
\begin{pf*}{Proof of Theorem \ref{free1}} Recalling (\ref
{Mar}), let $\tilde{\xi}_{ij}^{(n)}=(\xi_{ij}^{(n)}-\mu_n)/\sigma
_n$ for all $1\leq i< j\leq n$ and $n\geq2$. Then $\{\tilde{\xi
}_{ij}^{(n)}; 1\leq i < j \leq n, n\geq2\}$ satisfies (\ref
{shufic}) with $
\mu_n=0, \sigma_n=1$ and $p>4$. Let $\bd{\Delta}_{n,1}$ be
generated by $\{\tilde{\xi}_{ij}^{(n)}\}$ as in (\ref{Mar}). By
Lem\-ma~\ref{hahate}, almost surely,
%
%
\begin{equation}\label{boat}
F^{\bd{\Delta}_{n,1}/\sqrt{n}} \qquad\mbox{converges weakly to } \gamma_{M}
\end{equation}
as $n\to\infty$. It is easy to verify that
%
%
\begin{equation}\label{zhaod}
\bd{\Delta}_n=\underbrace{\sigma_n\bd{\Delta}_{n,1} +(n\mu_n)
\mathbf{I}_n}_{\bd{\Delta}_{n,2}} {}-{} \mu_n\mathbf{J}_n,
\end{equation}
where $\mathbf{I}_n$ is the $n\times n$ identity matrix, and $\mathbf{J}_n$
is the $n\times n$ matrix with all of its entries equal to 1.
Obviously, the eigenvalues of $\bd{\Delta}_{n,2}$ are $\sigma_n\cdot
\lambda_i(\bd{\Delta}_{n,1})+ n\mu_n,1\leq i \leq n$. By (\ref{boat}),
%
%
\begin{equation}\label{fangw}
\frac{1}{n}\sum_{i=1}^nI \biggl(\frac{\lambda_i(\bd{\Delta
}_{n,2})-n\mu_n}{\sqrt{n}\sigma_n}\leq x \biggr) \qquad\mbox{converges
weakly to } \gamma_{M}
\end{equation}
almost surely as $n\to\infty$. By (\ref{zhaod}) and the rank
inequality (see Lemma 2.2 from \cite{bai99}),
%
%
\begin{eqnarray}\label{qinmi}\qquad
&&\bigl\|F^{(\bd{\Delta}_n-n\mu_n\mathbf{I}_n)/\sqrt{n}\sigma_n}-F^{(\bd
{\Delta}_{n,2}-n\mu_n\mathbf{I}_n)/\sqrt{n}\sigma_n}\bigr\|
\nonumber\\[-8pt]\\[-8pt]
&&\qquad\leq
\frac{1}{n}\cdot\operatorname{rank} \biggl(\frac{\bd{\Delta}_n}{\sqrt
{n}\sigma_n}-\frac{\bd{\Delta}_{n,2}}{\sqrt{n}\sigma_n}
\biggr)
=\frac{1}{n}\cdot\operatorname{rank} \biggl(\frac{\mu
_n}{\sqrt{n}\sigma_n}\mathbf{J}_n \biggr)\leq\frac{1}{n}\to
0,\nonumber
\end{eqnarray}
where $\|f\|=\sup_{x\in\mathbb{R}}|f(x)|$ for any bounded,
measurable function $f(x)$ defined on $\mathbb{R}$. Finally, (\ref
{fangw}) and (\ref{qinmi}) lead to the desired conclusion.
\end{pf*}
\begin{pf*}{Proof of Theorem \ref{Becky}} Let
$\mathbf{V}_{n}=(v_{ij}^{(n)})$ be defined by
%
%
\begin{equation}\label{bymum}
v_{ii}^{(n)}=0 \quad\mbox{and}\quad v_{ij}^{(n)}=\frac{\xi_{ij}^{(n)}-\mu
_{n}}{\sigma_{n}}
\end{equation}
for any $1\leq i \ne j \leq n$ and $n\geq2$, where $\mathbf{A}_n=(\xi
_{ij}^{(n)})_{n\times n}$ as in (\ref{Dec}) with $\xi_{ii}^{(n)}=0$
for all $1\leq i \leq n$ and $n\geq2$.
It is easy to check that $\mathbf{A}_n=\mu_n(\mathbf{J}_n- \mathbf{I}_n) +
\sigma_n\mathbf{V}_n$, where all the entries of $\mathbf{J}_n$ are equal to
one, and $\mathbf{I}_n$ is the $n\times n$ identity matrix. Thus
\[
\frac{\mathbf{A}_n + \mu_n\mathbf{I}_n}{\sqrt{n}\sigma_n}-
\frac{\mathbf{V}_n}{\sqrt{n}}=
\frac{\mu_n\mathbf{J}_n}{\sqrt{n}\sigma_n}\qquad
\mbox{where all entries of $\mathbf{J}_n$ are equal to $1$.}
\]
By the rank inequality (see Lemma 2.2 from \cite{bai99}),
\[
\bigl\|F^{(\mathbf{A}_n+\mu_n\mathbf{I})/\sqrt{n}\sigma_n}-F^{n^{-1/2}
\mathbf{V}_{n}}\bigr\|\leq
\frac{1}{n}\cdot\operatorname{rank} \biggl(\frac{\mathbf{A}_n + \mu_n
\mathbf{I}_n}{\sqrt{n}\sigma_n}-\frac{\mathbf{V}_{n}}{\sqrt{n}} \biggr)\leq
\frac{1}{n}\to
0,
\]
where $\|f\|=\sup_{x\in\mathbb{R}}|f(x)|$ for any bounded,
measurable function $f(x)$ defined on $\mathbb{R}$ as in (\ref
{qinmi}). So, to prove the theorem, it is enough to show that
$F^{n^{-1/2}\mathbf{V}_{n}}$ converges weakly to
the semicircular law with the density given in statement of the
theorem. In view of normalization (\ref{bymum}), without loss of the
generality, we only need to prove the theorem under the conditions that
\[
E\omega_{ij}^{(n)}=0,\qquad E\bigl(\omega_{ij}^{(n)}\bigr)^{2}=1
\]
and
\[
\max
_{1\leq i<j\leq
n}E \bigl\{\bigl(\omega_{ij}^{(n)}\bigr)^{2}I \bigl(\bigl|\omega_{ij}^{(n)}\bigr|\geq
\varepsilon
\sqrt{n} \bigr) \bigr\}\to0 \qquad\mbox{as } n\to\infty
\]
for all $1\leq i, j\leq n$ and $n\geq2$. Given $\delta>0$, note that
\begin{eqnarray*}
&& \frac{1}{n^{2}\delta^{2}}
\sum_{1\leq i,j\leq n}
E \bigl\{\bigl(\omega_{ij}^{(n)}\bigr)^{2}I\bigl(\bigl|\omega_{ij}^{(n)}\bigr|\geq\delta
\sqrt{n}\bigr) \bigr\}\\
&&\qquad \leq \frac{2}{\delta^2}\cdot\max_{1\leq i<j\leq
n}E \bigl\{\bigl(\omega_{ij}^{(n)}\bigr)^{2}I \bigl(\bigl|\omega_{ij}^{(n)}\bigr|\geq
\delta
\sqrt{n} \bigr) \bigr\}\to0
\end{eqnarray*}
as $n\to\infty$. By Lemma \ref{missing} in the \hyperref[appendix]{Appendix}, $\tilde
{F}_n:=F^{n^{-1/2}\mathbf{V}_n}$, and hence $F^{n^{-1/2}(\mathbf{A}_{n}+\mu
_n\mathbf{I})}$,
converges weakly to the semicircular law.
\end{pf*}
\begin{pf*}{Proof of Corollary \ref{cheng}} To apply Theorem
\ref{Becky}, we first need to verify
%
%
\begin{equation}\label{niuren}
\max_{1\leq i<j\leq
n}E \bigl\{\bigl(\omega_{ij}^{(n)}\bigr)^{2}I \bigl(\bigl|\omega_{ij}^{(n)}\bigr|\geq
\varepsilon
\sqrt{n} \bigr) \bigr\}\to0
\end{equation}
as $n\to\infty$ for any $\varepsilon>0$, where $\omega
_{ij}^{(n)}:=(\xi_{ij}^{(n)}-\mu_n)/\sigma_n$. Note that $\mu
_n=p_n$ and $\sigma_n^2=p_n(1-p_n)$. Now, use the fact that $\xi
_{ij}^{(n)}$ take values one and zero only, and then the condition
$np_n(1-p_n)\to\infty$ to see that $|\omega_{ij}^{(n)}|\leq1/\sigma
_n=o (\sqrt{n} )$
as $n\to\infty$. Then (\ref{niuren}) follows. By Theorem \ref{Becky},
%
%
\begin{equation}\label{away}
\frac{1}{n}\sum_{i=1}^nI \biggl\{\frac{\lambda_i(\mathbf{A}_n) +
p_n}{\sqrt{np_n(1-p_n)}}\leq x \biggr\}
\end{equation}
converges weakly to the distribution with
density $\frac{1}{2\pi}\sqrt{4-x^{2}} I(|x|\leq2)$ almost surely. Notice
\[
\biggl\{\frac{\lambda_i(\mathbf{A}_n) + p_n}{\sqrt{np_n(1-p_n)}}\leq
x \biggr\}= \biggl\{\frac{\lambda_i(\mathbf{A}_n)}{\sqrt
{np_n(1-p_n)}}\leq x- \frac{p_n}{\sqrt{np_n(1-p_n)}} \biggr\}
\]
and $p_n/\sqrt{np_n(1-p_n)}\to0$ as $n\to\infty$. By using a
standard analysis, we obtain that, with probability one, $F^{
\mathbf{A}_n/\alpha_n}$ converges weakly to the semicircular law with density
$\frac{1}{2\pi}\sqrt{4-x^{2}} I(|x|\leq2)$, where $\alpha_n=\sqrt
{np_n(1-p_n)}$. Further, assume now $1/n \ll p_n\to0$ as $n\to\infty
$. Write
\[
\biggl\{\frac{\lambda_i(\mathbf{A}_n) + p_n}{\sqrt{np_n(1-p_n)}}\leq
x \biggr\}= \Biggl\{\frac{\lambda_i(\mathbf{A}_n)}{\sqrt{np_n}}\leq
x\sqrt{1-p_n}- \sqrt{\frac{p_n}{n}} \Biggr\}.
\]
Clearly, $ x\sqrt{1-p_n} - \sqrt{p_n/n}\to x$ as $n\to\infty$.
Thus, by (\ref{away}), we have $\frac{1}{n}\times\break\sum_{i=1}^{n}I\{\frac
{\lambda_i(\mathbf{A}_n)}{\sqrt{np_n}}\leq x\}$ converges weakly to the
semicircular law with density $\frac{1}{2\pi} 
\sqrt{4-x^{2}}
I(|x|\leq2)$.
\end{pf*}


We need the following lemma to prove Theorem \ref{Ruth}.
\begin{lemma}\label{whya} For $n\geq2$, let $\lambda_{n,1}\geq
\lambda
_{n,2}\geq\cdots\geq\lambda_{n,n}$ be real numbers. Set $\mu
_{n}=(1/n)\sum_{i=1}^n\delta_{\lambda_{n,i}}$. Suppose $\mu_n$
converges weakly to a probability measure $\mu$. Then, for any
sequence of integers $\{k_n; n\geq2\}$ satisfying $k_n=o(n)$ as
$n\to\infty$, we have
$\liminf_{n\to\infty}\lambda_{n,k_n}\geq\alpha$, where $\alpha
=\inf\{x\in\mathbb{R}\dvtx \mu([x, +\infty])=0\}$ with $\inf
\varnothing=+\infty$.
\end{lemma}
\begin{pf}
Since $\mu$ is a probability measure, we know that
$\alpha>-\infty$. Without loss of the generality, assume that $\alpha
>0$. For brevity of notation, write $k_n=k$. Set $\tilde{\mu
}_{n}=(n-k+1)^{-1}\sum_{i=k}^n\delta_{\lambda_{n,i}}$ for $n\geq k$.
Observe that
\[
\mu_n(B) -\tilde{\mu}_n(B)=\frac{1}{n}\sum_{i=1}^{k-1}I(\lambda
_{n,i}\in B) - \frac{k-1}{n(n-k+1)}\sum_{i=k}^{n}I(\lambda_{n,i}\in B)
\]
for any set $B\subset\mathbb{R}$, where $\sum_{i=1}^{k-1}I(\lambda
_{n,i}\in B)$ is understood to be zero if $k=1$. Thus, $|\mu_n(B)
-\tilde{\mu}_n(B)|\leq2k/n$.
Therefore,
%
%
\begin{equation}\label{ariza}
\tilde{\mu}_n \qquad\mbox{converges weakly to } \mu
\end{equation}
since $k=k_n=o(n)$ as $n\to\infty$. Easily,
\[
\lambda_{n,k}^mI(\lambda_{n,k}>0)\geq\frac{1}{n-k+1}\sum
_{i=k}^n\lambda_{n,i}^mI(\lambda_{n,i}>0)=\int_{0}^{\infty
}x^m\tilde{\mu}_n(dx)
\]
for any integer $m\geq1$. Write the last term above as $\int_{\mathbb
{R}}g(x)\tilde{\mu}_n(dx)$, where $g(x):=x^mI(x\geq0), x \in
\mathbb{R}$, is a continuous and nonnegative function. By (\ref
{ariza}) and the Fatou lemma,
%
%
\begin{equation}\label{hong1}\quad
\liminf_{n\to\infty}\lambda_{n,k}^mI(\lambda_{n,k}>0)\geq\liminf
_{n\to\infty}\int_{\mathbb{R}}g(x)\tilde{\mu}_n(dx)\geq\int
_{0}^{\infty}x^m\mu(dx) 
\end{equation}
for any $m\geq1$. If $\alpha<\infty$, then
%
%
\begin{equation}\label{hong2}
\int_{0}^{\infty}x^m\mu(dx)\geq\int_{\alpha-\varepsilon}^{\alpha
}x^m\mu(dx)\geq(\alpha-\varepsilon)^m\mu([\alpha-\varepsilon, \alpha])>0
\end{equation}
for any $\varepsilon\in(0, \alpha)$. Take the $(1/m)$th power for each
term in (\ref{hong1}) and (\ref{hong2}), and let $m\to\infty$ to get
\[
\liminf_{n\to\infty}\{\lambda_{n,k}I(\lambda_{n,k}>0)\} \geq
\alpha- \varepsilon
\]
for any $\varepsilon\in(0, \alpha)$. By sending $\varepsilon\downarrow0$
and using the fact $\alpha>0$, the conclusion is yielded.

If $\alpha=+\infty$, notice
\[
\int_{0}^{\infty}x^m\mu(dx)\geq\int_{\rho}^{\infty}x^m\mu
(dx)\geq\rho^m\mu([\rho, \infty))>0
\]
for any $\rho>0$. Using the same argument as above and then letting
$\rho\to+\infty$, we get the desired assertion.
\end{pf}
\begin{pf*}{Proof of Lemma \ref{Mikea}}
(i) By Theorem \ref{Becky}, $F^{n^{-1/2}\mathbf{U}_n}$
converges weakly to the semicircular law with density function $\frac
{1}{2\pi}\sqrt{4-x^{2}}I(|x|\leq2)$. Use Lemma \ref{whya} to have that
%
%
\begin{equation}
\liminf_{n\to\infty}\frac{\lambda_{\max}(\mathbf{U}_n)}{\sqrt{n}}
\geq2\qquad\mbox{a.s.}
\end{equation}
Now we prove the upper bound, that is,
%
%
\begin{equation}\label{upper2}
\limsup_{n\to\infty}\frac{\lambda_{\max}(\mathbf{U}_n)}{\sqrt{n}}\leq2
\qquad\mbox{a.s.}
\end{equation}
Define
\[
\delta_{n}=\frac{1}{\log(n+1)},\qquad
\tilde{u}_{ij}^{(n)}=u_{ij}^{(n)}I\bigl(\bigl|u_{ij}^{(n)}\bigr|\leq
\delta_{n}\sqrt{n}\bigr)\quad \mbox{and} \quad\tilde{\mathbf{U}}_{n}=\bigl(\tilde
{u}_{ij}^{(n)}\bigr)_{1\leq i, j\leq n}
\]
for $1\leq i\leq j\leq n$ and $n\geq1$. By the Markov inequality,
\begin{eqnarray*}
P(\mathbf{U}_n\neq\tilde{\mathbf{U}}_{n})&\leq&P\bigl(\bigl|u_{ij}^{(n)}\bigr| >
\delta_{n}\sqrt{n} \mbox{ for some } 1\leq i,j\leq n\bigr) \\
&\leq&n^{2}\max_{1\leq i,j\leq n}P\bigl(\bigl|u_{ij}^{(n)}\bigr|>
\delta_{n}\sqrt{n}\bigr)\\
&\leq&\frac{K
(\log(n+1))^{6+\delta}}{n^{1+(\delta/2)}},
\end{eqnarray*}
where $K=\sup_{1\leq i,j\leq
n, n\geq1}E|u_{ij}^{(n)}|^{6+\delta}<\infty$. Therefore, by the
Borel--Cantelli lemma,
%
%
\begin{equation}\label{JiL}
P(\mathbf{U}_n=\tilde{\mathbf{U}}_{n} \mbox{ for sufficiently large } n)=1.
\end{equation}
From $Eu_{ij}^{(n)}=0$, we have that
%
%
\begin{equation}\label{hood}\qquad
\bigl|Eu_{ij}^{(n)}I\bigl(\bigl|u_{ij}^{(n)}\bigr|\leq
\delta_{n}\sqrt{n}\bigr)\bigr|=\bigl|Eu_{ij}^{(n)}I\bigl(\bigl|u_{ij}^{(n)}\bigr|>
\delta_{n}\sqrt{n}\bigr)\bigr| \leq\frac{K}{
(\delta_{n}\sqrt{n})^{5+\delta}}
\end{equation}
for any $1\leq i\leq j\leq n, n\geq1$. Note that $\lambda_{\max}(
\mathbf{A}+\mathbf{B})\leq\lambda_{\max}(\mathbf{A}) + \lambda_{\max}(\mathbf{B})$,
and $\lambda_{\max}(\mathbf{A})\leq\|\mathbf{A}\|\leq n\cdot\max_{1\leq i,
j\leq
n}|a_{ij}|$ for any $n\times n$ symmetric matrices $\mathbf{A}=(a_{ij})$
and $\mathbf{B}$. We have from (\ref{hood}) that
\begin{eqnarray*}
&&\lambda_{\max}(\tilde{\mathbf{U}}_{n})-\lambda_{\max}\bigl(\tilde{
\mathbf{U}}_{n}-E(\tilde{\mathbf{U}}_{n})\bigr)\\
&&\qquad\leq\lambda_{\max}(E\tilde{
\mathbf{U}}_{n})
\leq n\max_{1\leq i, j\leq n}\bigl|Eu_{ij}^{(n)}I\bigl(u_{ij}^{(n)}
\leq
\delta_{n}\sqrt{n}\bigr)\bigr|\\
&&\qquad\leq\frac{K}{\delta_{n}^{5+\delta}(\sqrt
{n})^{3+\delta}}
\end{eqnarray*}
for any $n\geq1$. This and (\ref{JiL}) imply that
\[
\limsup_{n\to\infty}\frac{\lambda_{\max}(\mathbf{U}_n)}{\sqrt
{n}}=\limsup_{n\to\infty}
\frac{\lambda_{\max}(\tilde{\mathbf{U}}_n)}{\sqrt{n}} \leq
\limsup_{n\to\infty}\frac{\lambda_{\max}(\tilde{
\mathbf{U}}_{n}-E\tilde{\mathbf{U}}_{n})}{\sqrt{n}}
\]
almost surely.\vspace*{1pt} Note that $|\tilde{u}_{ij}^{(n)}|\leq|u_{ij}^{(n)}|$
and $\operatorname{Var}(\tilde{u}_{ij}^{(n)})\leq E(u_{ij}^{(n)})^{2}=1$, to
save notation, without loss of generality, we will prove (\ref
{upper2}) by assuming that
\[
E\bigl(u_{ij}^{(n)}\bigr)=0,\qquad
E\bigl(u_{ij}^{(n)}\bigr)^{2}\leq1,\qquad \bigl|u_{ij}^{(n)}\bigr|\leq\frac{2\sqrt{n}}{\log
(n+1)}
\]
and
\[
\max_{1\leq i,j\leq n,n\geq1} E\bigl|u_{ij}^{(n)}\bigr|^{6+\delta
}<\infty
\]
for all $1\leq i, j\leq n$ and $n\geq1$. Now,
\[
\max_{i,j,n}
E\bigl|u_{ij}^{(n)}\bigr|^{3} \leq\max_{i,j,n}\bigl(
E\bigl|u_{ij}^{(n)}\bigr|^{6+\delta}\bigr)^{3/(6+\delta)}=K^{3/(6+\delta)}
\]
by the H\"{o}lder inequality. Hence,
%
%
\begin{equation}
\max_{1\leq i, j\leq n}E\bigl|u_{ij}^{(n)}\bigr|^{l}\leq K^{3/(6+\delta)}\cdot
\biggl(\frac{2\sqrt{n}}{\log(n+1)} \biggr)^{l-3}
\end{equation}
for all $n\geq1$ and $l\geq3$, where $K$ is a constant. The
inequality in (\ref{upper2}) follows from Lemma \ref{sprc} in the
\hyperref[appendix]{Appendix}. Thus the first limit in the lemma is proved. Applying this
result to $-\mathbf{U}_n$, we obtain
%
%
\begin{equation}
\lim_{n\to\infty}\frac{\lambda_{\min}(\mathbf{U}_n)}{\sqrt{n}}
=-\lim_{n\to\infty}\frac{\lambda_{\max}(-\mathbf{U}_n)}{\sqrt{n}}=-2
\qquad\mbox{a.s.}
\end{equation}
Since $\|\mathbf{U}_n\|=\max\{\lambda_{\max}(\mathbf{U}_n), -\lambda
_{\min}(\mathbf{U}_n)\}$, the above and the first limit in the lemma yield
the second limit.

(ii) Let $\hat{\mathbf{U}}_n=\mathbf{U}_n-\operatorname{diag}(u_{ii}^{(n)})_{1\leq
i \leq n}$. It is not difficult to check that both $|\lambda
_{\max}(\hat{\mathbf{U}}_n)-\lambda_{\max}(\mathbf{U}_n)|$ and $| \|\hat
{\mathbf{U}}_n\|-\|\mathbf{U}_n\| |$ are bounded by $\|
{\operatorname{diag}}(u_{ii}^{(n)})_{1\leq i \leq n}\|=\max_{1\leq i \leq
n}|u_{ii}^{(n)}|$. By (i), it is enough to show
%
%
\begin{equation}\label{snowin}
\max_{1\leq i \leq n}\bigl|u_{ii}^{(n)}\bigr|/n^{1/3}\to0 \qquad\mbox{a.s.}
\end{equation}
as $n\to\infty$. In fact, by the Markov inequality
\begin{eqnarray*}
\sum_{n=1}^{\infty}P\Bigl(\max_{1\leq i \leq n}\bigl|u_{ii}^{(n)}\bigr| \geq
n^{1/3}t\Bigr) & \leq& \sum_{n=1}^{\infty}n\cdot\max_{1\leq i \leq
n}P\bigl(\bigl|u_{ii}^{(n)}\bigr| \geq n^{1/3}t\bigr)\\
& \leq& \sum_{n=1}^{\infty}\frac{t^{-6-\delta}}{n^{1+(\delta
/3)}}\cdot\sup_{1\leq i,j\leq
n, n\geq1}E\bigl|u_{ij}^{(n)}\bigr|^{6+\delta}<\infty
\end{eqnarray*}
for any $t>0$. Thus, (\ref{snowin}) is concluded by the
Borel--Cantelli lemma.
\end{pf*}
\begin{pf*}{Proof of Theorem \ref{Ruth}} Let $\mathbf{J}_{n}$ be
the $n\times n$ matrix whose $n^2$
entries are all equal to 1. Let $\mathbf{V}_n$ be defined as in (\ref
{bymum}). Then
$\mathbf{B}_n:=\mathbf{A}_n+\mu_n\mathbf{I}_n=\sigma_{n}\mathbf{V}_n+\mu_{n}
\mathbf{J}_{n}$. First, by Lemma \ref{Mikea},
%
%
\begin{equation}\label{erya}
\lim_{n\to\infty}\frac{\lambda_{\max}(\mathbf{V}_n)}{\sqrt{n}}=2
\qquad\mbox{a.s.}\quad\mbox{and}\quad
\lim_{n\to\infty}\frac{\|\mathbf{V}_n\|}{\sqrt{n}}=2\qquad
\mbox{a.s.}
\end{equation}
Since $\mathbf{V}_n$ is symmetric, $\|\mathbf{V}_n\|={\sup_{x\in\mathbb
{R}^n\dvtx\|x\|=1}}\|\mathbf{V}_nx\|={\sup_{\|x\|=1}}|x^{T}\mathbf{V}_nx|$.
By definition
%
%
\begin{eqnarray}\label{laos}
\lambda_{\max}(\mathbf{B}_n)&=&\sup_{\|x\|=1}\{\sigma_{n}(x^{T}
\mathbf{V}_nx)+\mu_{n}(x^{T}\mathbf{J}_{n}x)\}\nonumber\\[-8pt]\\[-8pt]
&=&\sup_{\|x\|=1}\{\sigma_{n}(x^{T}\mathbf{V}_nx)+\mu_{n}(
\mathbf{1}'x)^{2}\},\nonumber
\end{eqnarray}
because $\mathbf{J}=\mathbf{1}\cdot\mathbf{1}^T$, where $\mathbf{1}=(1,\ldots
,1)^T\in\mathbb{R}^{n}$. Second, by Theorem \ref{Becky},
$F^{\mathbf{B}_n/\sqrt{n}\sigma_n}$ converges weakly to the semicircular law
$\frac{1}{2\pi}\sqrt{4-x^{2}}I(|x|\leq2)$. From Lemma \ref{whya},
we know that
%
%
\begin{equation}\label{hahaa}
\liminf_{n\to\infty}\frac{\lambda_{k_n}(\mathbf{B}_n)}{\sqrt
{n}\sigma_n} \geq2\qquad\mbox{a.s.}
\end{equation}

Now we are ready to prove the conclusions.

\begin{longlist}
\item
It is easy to check that
$\sup_{\|x\|=1}\{(\mathbf{1}^{\prime}x)^{2}\}=n$. By (\ref{laos}), $\lambda
_{\max}(\mathbf{B}_n)\leq\sigma_n\|\mathbf{V}_n\| + n|\mu_n|$. Thus
$\limsup_{n\to\infty}\lambda_{\max}(\mathbf{B}_n)/\sqrt{n}\sigma
_n\leq2$ a.s. by (\ref{erya}) under the assumption $\mu
_n/(n^{-1/2}\sigma_n)\to0$ as $n\to\infty$. Since $\lambda
_{\max}(\mathbf{B}_n)=\mu_n +\lambda_{\max}(\mathbf{A}_n)$. From (\ref
{hahaa}) we see that $\lim_{n\to\infty}\lambda_{k_n}(
\mathbf{A}_n)/\sqrt{n}\sigma_n = 2$ a.s. when $\mu_n/(n^{-1/2}\sigma
_n)\to0$ as $n\to\infty$. In particular, $\lim_{n\to\infty
}\lambda_{\max}(\mathbf{A}_n)/\sqrt{n}\sigma_n = 2$ a.s. Under the same
condition, we also have $\lim_{n\to\infty}\lambda_{\max}(-
\mathbf{A}_n)/\sqrt{n}\sigma_n = 2$ a.s. Finally, using $\|\mathbf{A}_n\|=\max
\{\lambda_{\max}(\mathbf{A}_n), \lambda_{\max}(-\mathbf{A}_n)\}$, we obtain
that $\lim_{n\to\infty}\|\mathbf{A}_n\|/\sqrt{n}\sigma_n = 2$ a.s.

\item Without loss of generality, assume $\mu_{n}>0$ for all $n\geq2$.
From (\ref{laos}) we see that
\begin{eqnarray*}
&&\mu_{n}\sup_{\|x\|=1}\{(\mathbf{1}^{\prime}x)^{2}\}-\sigma_{n}\sup_{\|x\|
=1}\{|x^{T}\mathbf{V}_nx|\}
\\
&&\qquad\leq
\lambda_{\max}(\mathbf{B}_n)
\leq\mu_{n}\sup_{\|x\|=1}\{(\mathbf{1}^{\prime}x)^{2}\}+
\sigma_{n}\sup_{\|x\|=1}\{|x^{T}\mathbf{V}_nx|\}.
\end{eqnarray*}
Hence, $n\mu_n - \sigma_n\|\mathbf{V}_n\|\leq\lambda_{\max}(
\mathbf{B}_n)\leq n\mu_n+\sigma_n\|\mathbf{V}_n\|$. Consequently, if $\mu
_{n}\gg
n^{-1/2}\sigma_{n}$, by (\ref{erya}), we have
\[
\lim_{n\to\infty}\frac{\lambda_{\max}(\mathbf{A}_n)}{n\mu_{n}}=
\lim_{n\to\infty}\frac{\lambda_{\max}(\mathbf{B}_n)}{n\mu_{n}}=1\qquad
\mbox{a.s.}
\]
since $\lambda_{\max}(\mathbf{B}_n)=\mu_n +\lambda_{\max}(\mathbf{A}_n)$.

\item Since $\mathbf{B}_n=\sigma_{n}\mathbf{V}_n+\mu_{n}\mathbf{J}_{n}$ and $\|
\mathbf{J}_n\|=n$, by the triangle inequality of $\|\cdot\|$,
\[
n|\mu_n|-\sigma_n\|\mathbf{V}_n\| \leq\|\mathbf{B}_n\| \leq n|\mu
_n|+\sigma_n\|\mathbf{V}_n\|.
\]
By (\ref{erya}) and the definition that $\mathbf{A}_n=\mathbf{B}_n-\mu_n\mathbf{I}_n$, we obtain
\[
\lim_{n\to\infty}\frac{\|\mathbf{A}_n\|}{n|\mu_n|}
=\lim_{n\to\infty}\frac{\|\mathbf{B}_n\|}{n|\mu_n|}=1\qquad\mbox{a.s.}
\]
as $|\mu_n|\gg n^{-1/2}\sigma_n$.\qed
\end{longlist}
\noqed\end{pf*}\vspace*{-14pt}

\begin{appendix}\label{appendix}
\section*{Appendix}

\setcounter{lemma}{0}
\begin{lemma}[(Sakhanenko)]\label{mamam} Let $\{\xi_{i};i=1,2,\ldots\}$
be a sequence of
independent random variables with mean zero and variance
$\sigma_{i}^{2}$. If $E|\xi_{i}|^{p}<\infty$ for some $p>2$, then
there exists a constant $C>0$ and $\{\eta_{i};i=1,2,\ldots\}$, a
sequence of independent normally distributed random variables with
$\eta_{i}\sim N(0,\sigma_{i}^{2})$ such that
\[
P\Bigl({\max_{1\leq k\leq n}}|S_{k}-T_{k}|>x\Bigr)\leq
\frac{C}{1+|x|^{p}}\sum_{i=1}^{n}E|\xi_{i}|^{p}
\]
for any n and $x>0$, where $S_{k}=\sum_{i=1}^{k}\xi_{i}$ and
$T_{k}=\sum_{i=1}^{k}\eta_{i}$.
\end{lemma}

Let $\mathbf{W}_n=(\omega_{ij}^{n})_{1\leq i,j\leq n}$ be an $n\times n$
symmetric matrix, where $\{\omega_{ij}^{n}; 1\leq i\leq j\leq n\}$
are random variables defined on the same probability space. We need the
following two results from Bai \cite{bai99}.
\begin{lemma}[(Theorem 2.4 in \cite{bai99})]\label{missing}
For each $n\geq2$, let $\{\omega_{ij}^{n}; 1\leq i\leq j\leq
n\}$ be independent random variables (not necessarily
identically distributed) with $\omega_{ii}^{n}=0$ for all $1\leq i
\leq n$, $E(\omega_{ij}^{n})=0$ and
$E(\omega_{ij}^{n})^{2}=\sigma^{2}>0$ for all $1\leq i< j\leq n$,
and
\[
\lim_{n\to\infty}\frac{1}{n^{2}\delta^{2}}\sum_{1\leq
i,j\leq n}E(\omega_{ij}^{n})^{2}I\bigl(|\omega_{ij}^{n}|\geq
\delta\sqrt{n}\bigr)=0
\]
for any $\delta>0$. Then $F^{n^{-1/2}\mathbf{W}_n}$ converges weakly to
the semicircular law
of scale-parameter $\sigma$ with density function
%
%
\setcounter{equation}{0}
\begin{equation}
p_{\sigma}(x) = \cases{
\dfrac{1}{2\pi\sigma^{2}}\sqrt{4\sigma^{2}-x^{2}},&\quad if $|x|\leq
2\sigma$;\vspace*{2pt}\cr
0,&\quad otherwise.}
\end{equation}
\end{lemma}

Some recent results in \cite{BZ2008,PGZ2007} are in the realm of the
above lemma.
\begin{lemma}[(Remark 2.7 in \cite{bai99})]\label{sprc}
Suppose, for each $n\geq1$, $\{\omega_{ij}^{(n)}; 1\leq i\leq j\leq
n\}$ are independent random variables (not necessarily identically
distributed) with mean $\mu=0$ and variance no larger than
$\sigma^{2}$. Assume there exist constants $b>0$ and
$\delta_{n}\downarrow0$ such that $\sup_{1\leq i,j\leq
n}E|\omega_{ij}^{(n)}|^{l}\leq b(\delta_{n}\sqrt{n})^{l-3}$ for all
$n\geq1$ and $l\geq3$. Then
\[
\limsup_{n\to\infty}\frac{\lambda_{\max}(\mathbf{W}_n)}{n^{1/2}}\leq
2\sigma\qquad\mbox{a.s.}
\]
\end{lemma}
\end{appendix}


%
\printaddresses

\end{document}